\documentclass [11pt,reqno]{amsart}
\usepackage{amscd,amsmath,amssymb,euscript,graphics}
\usepackage[frame,cmtip,curve,arrow,matrix,line,graph]{xy}
\usepackage{hyperref}

\title{Stability conditions and Stokes factors}
\author{Tom Bridgeland}
\address{
All Souls College,
Oxford,
UK}
\email{bridgeland@maths.ox.ac.uk}
\thanks{T.B. supported by a Royal Society University Research Fellowship}
\author{Valerio Toledano Laredo}
\address{Department of Mathematics,
Northeastern University,
567 Lake Hall,
360 Huntington Avenue,
Boston MA 02115.
\hfill\break\indent
School of Mathematics,
Institute for Advanced Study,
Princeton, NJ 08540.}
\email{V.ToledanoLaredo@neu.edu}
\thanks{V.T.L. supported in part by NSF grants DMS--0707212 and DMS--0635607}
\dedicatory{Dedicated to Graeme Segal, whose approach to mathematics has been an inspiration to both of us}

\newtheorem*{thm}{Theorem}

\newtheorem*{cor}{Corollary}
\newtheorem*{prop}{Proposition}
\newtheorem*{lemma}{Lemma}

\newenvironment{pf}{\paragraph{{\sc Proof}}}{\qed\par\medskip}
\theoremstyle{definition}
\newtheorem*{defn}{Definition}
\newtheorem*{remark}{Remark}

\newtheorem*{example}{Example}

\newcommand {\bb}{\mathfrak{b}}
\newcommand{\g}{\mathfrak{g}}
\newcommand{\h}{\mathfrak{h}}

\newcommand{\n}{\mathfrak{n}}

\newcommand{\GL}{\operatorname{GL}}
\renewcommand{\leq}{\leqslant}
\renewcommand{\geq}{\geqslant}
\newcommand{\ad}{\operatorname{ad}}

\newcommand{\End}{\operatorname{End}}
\newcommand{\D}{\operatorname{\mathcal{D}}}

\newcommand{\isom}{\cong}

\newcommand{\C}{\mathbb C}

\newcommand{\Z}{\mathbb Z}
\newcommand{\A}{\mathcal A}

\newcommand{\Flag}{\operatorname{Flag}}
\newcommand{\Hom}{\operatorname{Hom}}

\newcommand{\lra}{\longrightarrow}
\newcommand{\tensor}{\otimes}
\newcommand{\R}{\mathbb{R}}

\newcommand{\reg}{{\operatorname{reg}}}
\newcommand{\M}{\mathcal{M}}
\renewcommand{\H}{{\mathcal H}}

\newcommand{\hreg}{\h_{\operatorname{reg}}}

\newcommand{\Stab}{\operatorname{Stab}}
\renewcommand{\Im}{\operatorname{Im}}

\renewcommand{\Re}{\operatorname{Re}}

\newcommand{\Rep}{\operatorname{Rep}}
\newcommand{\god}{\g_{\text{od}}}

\newcommand{\gl}{\mathfrak{gl}}

\newcommand{\wh}[1]{\widehat{#1}}

\newcommand{\IC}{\mathbb C}
\newcommand{\IH}{\mathbb H}

\newcommand{\IN}{\mathbb N}
\newcommand{\IP}{\mathbb P}
\newcommand{\IQ}{\mathbb Q}
\newcommand{\IR}{\mathbb R}

\newcommand{\IZ}{\mathbb Z}

\renewcommand {\L}{\mathcal L}
\newcommand {\calR}{\mathcal R}
\newcommand {\calS}{\mathcal S}
\newcommand {\U}{\mathcal U}

\newcommand{\Ker}{\operatorname{Ker}}

\newcommand{\eg}{{\it e.g. }}

\newcommand{\ol}[1]{\overline{#1}}

\newcommand{\wrt}{with respect to }
\newcommand{\aand}{\qquad\text{and}\qquad}

\newcommand {\bfX}{\mathbb{X}}
\newcommand {\rep}{representation }
\newcommand {\clockwise}{\stackrel{\curvearrowright}{\prod}}
\newcommand {\Mod}{\operatorname{Mod}}

\newcommand {\PhiZ}{\Phi} %to be reverted to {\Phi^Z} if we decide to assume that Z is not regular
\newcommand {\LambdaZ}{\Lambda} %to be reverted to {\Lambda^Z} if we decide to assume that Z is not regular
\renewcommand {\Vec}{\operatorname{Vec}}%{\mathbf{Vec}}%
\newcommand {\RRep}{\operatorname{Rep}}%{\mathbf{Rep}}

\newcommand {\BA}{\wh{B}(\A)}
\newcommand {\CA}{\wh{C}(\A)}
\newcommand {\HA}{H(\A)}
\newcommand {\KA}{K(\A)}
\newcommand {\NA}{\wh{N}(\A)}
\newcommand {\hA}{\h(\A)}
\newcommand {\nA}{\wh{\n}(\A)}
\newcommand {\CbarA}[1]{C_{\leq #1}(\A)}%{\overline{C}_{\leq #1}(\A)}
\newcommand {\nbarA}[1]{\n_{\leq #1}(\A)}%{\overline{\n}_{\leq #1}(\A)}
\newcommand {\NbarA}[1]{N_{\leq #1}(\A)}%{\overline{\n}_{\leq #1}(\A)}
\newcommand {\bbarA}[1]{\bb_{\leq #1}(\A)}
\newcommand {\BbarA}[1]{B_{\leq #1}(\A)}
\newcommand {\PbarA}[1]{P_{\leq #1}}%(\A)}

\newcommand{\fd}{finite--dimensional }

\newcommand{\RH}{Ringel--Hall }

\newcommand {\HN}{Harder--Narasimhan }

\renewcommand{\SS}{\operatorname{SS}}
\newcommand {\Splus}{\operatorname{1_{\A}}}

\newcommand {\dchi}{d\chi}

\begin{document}

\begin{abstract}
Let $\A$ be the category of modules over a complex, \fd algebra.
We show that the space of stability conditions on $\A$ parametrises
an isomonodromic family of irregular connections on $\IP^1$ with
values in the Hall algebra of $\A$. The residues of these connections
are given by the holomorphic generating function for counting
invariants in $\A$ constructed by D. Joyce \cite{Joyce}.
\end{abstract}
\maketitle

\setcounter{tocdepth}{1}
\tableofcontents

\section{Introduction}
%===============

\subsection{}
%--------------

This paper stems from our attempt to understand the recent work of Joyce
\cite{Joyce} on holomorphic generating functions for counting invariants in
an abelian category $\A$. Somewhat unexpectedly, a more conceptual
understanding of Joyce's formulae can be obtained by viewing them as
defining  an irregular connection on $\IP^1$ with values in the Ringel--Hall
Lie algebra of $\A$. This leads to a picture whereby stability conditions on
$\A$ can be naturally interpreted as defining Stokes data for such connections.

We begin with a leisurely introduction reviewing the salient points
of Joyce's work and summarising our main results, a more precise
formulation of which may be found in the body of the paper. From
\ref{ss:Stokes from now} on, our exposition assumes a passing
knowledge of Stokes phenomena; the reader can find an introduction
to this material in Section \ref{se:Irregular connections}.

\subsection{}
%--------------

It is a familiar fact from Geometric Invariant Theory that if one wants to form moduli
spaces parameterising algebro--geometric objects such as coherent sheaves or
modules over an algebra, one first needs to restrict to some subclass of (semi)stable
ones. The required notion of stability is usually not given {\it a priori},
%(Mumford stability for bundles on curves is perhaps an exception in this respect)
but rather corresponds
to some particular choice of weights. As these weights  vary, the corresponding
subclasses of semistable objects undergo discontinuous changes; in many cases
the space of all possible weights has a wall--and--chamber decomposition such
that the subclass of semistable objects is constant in each chamber but jumps
as one moves across a wall.

More recently, these spaces of weights, or stability conditions, have been studied as
interesting objects in their own right. Spaces of stability conditions on triangulated
categories were introduced by the first author in \cite{Bridgeland} following earlier
work by M. Douglas \cite{Douglas}. Considerations from Mirror Symmetry suggest
that these spaces should have interesting geometric structures closely related to
Frobenius structures \cite{Bridgeland2}; in particular one expects that they should
carry natural families of irregular connections. Recently, progress has been made
towards defining such structures \cite{Joyce, KS} although the picture is still far from
clear.

In this paper we shall be concerned with stability conditions on
an abelian category $\A$, and in fact with the special case when
$\A$ is the category of \fd modules over a \fd algebra.  We hope
that, with further study these ideas can be applied to more
 general situations such as spaces of stability conditions on
derived categories of coherent sheaves.
%; as we explain below
%however, such an extension will require several new ideas.

\subsection{}
%--------------

Suppose then that $\A$ is the abelian category of \fd modules over
a \fd associative $\C$--algebra $R$. Let $K(\A)$ be the Grothendieck
group of $\A$ and $K_{> 0}(\A)\subset K(\A)$ the positive cone
spanned by the classes of nonzero modules. Let $\IH\subset \C$
denote the upper half--plane.  For our purposes, a {\it stability
condition} on $\A$ is just a homomorphism of abelian groups
$Z\colon K(\A)\to \C$ such that
$$Z(K_{>0}(\A))\subset\IH.$$
In other words, a stability condition is a choice $Z(M)\in\IH$ for each nonzero module
$M$ such that $Z$ is additive on short exact sequences. Given such a stability condition
$Z$ each nonzero module $M$ has a well-defined {\it phase}
$$\phi(M)=\frac{1}{\pi}\arg Z(M) \in (0,1),$$
and a non--zero module $M$ is said to be {\it semistable} with respect to $Z$ if every
non--zero submodule $A\subset M$ satisfies $\phi(A)\leq\phi(M)$.

Since the category $\A$ has finite length, the Grothendieck group $K(\A)\isom
\Z^N$ is freely generated by the classes of the simple modules, and
the space of all stability conditions $\Stab(\A)$ can be identified with the complex
manifold $\IH^N$. It is easy to see that for each class $\alpha\in K_{> 0}(\A)$,
there is a finite collection of codimension--one real submanifolds of $\Stab(\A)$
such that in each connected component of their complement the set of semistable
modules of class $\alpha$ is constant. This is the wall--and--chamber structure
referred to above. Moreover, it follows from results of King \cite{King} that for
any stability condition $Z$ and class $\alpha\in K_{>0}(\A)$, there is a
projective scheme which is a coarse moduli space for semistable modules
of type $\alpha$. We can thus view this simple--minded example as a good
model for studying wall--crossing phenomena.

\subsection{}\label{ss:Hall}
%--------------

From an algebraic perspective, these wall--crossing phenomena give rise to,
and may be studied as change of bases within the {\it Ringel--Hall algebra}
$\H(\A)$ of $\A$, an idea which has its origins in the work of Reineke
\cite{Reineke}, and was greatly developed by Joyce \cite{Joyce1,Joyce4}.
For a survey of Hall algebras over finite fields see \cite{Schiff}. The variant
we shall use here was sketched by Kapranov and Vasserot \cite{KV} and
described in detail by Joyce \cite{Joyce1}.

Consider first the vector space $\H(\A)$ of complex--valued
constructible functions on the moduli stack of all $R$--modules.
This vector space can be endowed with an associative product $*$
for which
$$(f_1*\cdots *f_n)(M)=\int f_1(M_1/M_0) \cdots f_n(M_n/M_{n-1})\dchi,$$
where the Euler characteristic integral is taken over the variety %$\Flag^n(M)$
parameterising flags
$$0=M_0\subset M_1\subset \cdots \subset M_n=M$$
of submodules of $M$ of length $n$. The characteristic function $1_0$
of the zero module is the identity element. The resulting algebra has an
obvious grading
\[\H(\A)=\bigoplus_{\alpha\in K_{\geq 0}(\A)}\H_\alpha(\A)\]
where $\H_\alpha(\A)$ is the space of functions supported on
modules of class $\alpha$.

The algebra $\H(\A)$ is usually much too big, and one considers
instead the subalgebra $C(\A)$ generated by the characteristic
functions $\kappa_\alpha$ of the sets of modules of class
$\alpha$ as $\alpha$ varies in $K_{> 0}(\A)$. The rule
$$\Delta(f)(M,N)=f(M\oplus N)$$ %,\quad f\in\H(\A)$$
defines a coproduct $\Delta:C(\A)\to C(\A)\otimes C(\A)$ on $C(\A)$
which endows it with the structure of a cocommutative bialgebra.
The corresponding Lie algebra of primitive elements $\n(\A)$
coincides with the space of functions supported on indecomposable
objects, and the inclusion $\n(\A)\subset C(\A)$ identifies $C(\A)$
with the universal enveloping algebra of $\n(\A)$.

The grading on $\H(\A)$ induces gradings on $C(\A)$ and $\n(\A)$.
One can use the grading on $\n(\A)$ to define an extended Lie
algebra $\bb(\A)=\h(\A)\ltimes\n(\A)$ by endowing
$$\h(\A)=\Hom_{\Z}(K(\A),\C)$$
with a trivial bracket, and setting $[Z,f]=Z(\alpha)f$ for $Z\in \h(\A)$
and $f\in \n_\alpha(\A)$.

It will also be necessary in what follows to consider the completions
$\nA$ and $\wh{\bb}(\A)$ of $\n(\A)$ and $\bb(\A)$ with respect to
their gradings. We shall collectively refer to these Lie algebras as
\RH Lie algebras of $\A$.

\subsection{}\label{ss:Joyce review}
%--------------

The original motivation of the present work was to understand Joyce's
remarkable paper \cite{Joyce}, some features of which we now briefly
recall.

Given a stability condition $Z\in\Stab(\A)$, the characteristic function of
the set of semistable modules of a given class $\gamma\in K_{>0}(\A)$
defines an element $\delta_\gamma\in C_\gamma(\A)$ which plainly
encodes the discontinuous behaviour of that set. Joyce defines closely
related elements $\epsilon_\alpha$ by the (finite) sum
\[\epsilon_\alpha=
\sum_{n\geq 1}
\sum_{\substack{\gamma_1+\cdots+\gamma_n=\alpha\\ Z(\gamma_i)\in\R_{>0}Z(\alpha)}}
\frac{(-1)^{n-1}}{n}\,\delta_{\gamma_1}*\cdots*\delta_{\gamma_n},\]
and proves that they lie in the Lie algebra $\n(\A)$.
%,and that the original functions $\delta_\gamma$ may
%be recovered as
%$$\delta_\gamma=
%\sum_{n\geq 1}
%\sum_{\substack{\alpha_1+\cdots+\alpha_n=\gamma\\ Z(\alpha_i)\in\R_{>0}Z(\gamma)}}
%\frac{\epsilon_{\alpha_1}*\cdots*\epsilon_{\alpha_n}}{n!}.$$
%Note that the above sums are finite because the elements $\alpha$ and
%$\gamma$ lie in the positive cone $K_{>0}(\A)\cong\IN^N\setminus\{0\}$.}
Considered as a function $\Stab(\A)\to\n_\alpha(\A)$, $\epsilon_\alpha$ is
constant on the chambers in $\Stab(\A)$ defined by $\alpha$ and exhibits
discontinuous behaviour on their walls.

Joyce then considers, for any $\alpha\in K_{>0}(\A)$, a generating function
$f_\alpha\colon \Stab(\A)\to \n_\alpha(\A)$ given by a Lie series of the form
\begin{equation}
\label{liedseries}
f_\alpha=
\sum_{n\geq 1}
\sum_{\substack{\alpha_i\in K_{>0}(\A)\\ \alpha_1+\cdots+\alpha_n=\alpha}}
F_n(Z(\alpha_1),\ldots, Z(\alpha_n))\,
\epsilon_{\alpha_1}*\cdots *\epsilon_{\alpha_n},
\end{equation}
where $F_n\colon (\C^*)^n\to \C$ is a function on $n$ complex variables,
with $F_1\equiv 1$. He then proves the following result (see \cite[\S 3]
{Joyce} for a more precise statement of (i)--(ii)).

\begin{thm}[Joyce]\hfill
\begin{enumerate}
\item The functions $F_n$ can be chosen to be holomorphic with branchcuts
which precisely balance the discontinuities of the $\epsilon_\alpha$, thus
resulting in a continuous, holomorphic function $f_\alpha$, independently
of which algebra $R$ (and hence which abelian category $\A$) one starts
with.
\item The functions $F_n$ are uniquely characterised by the above
requirements and a few additional mild assumptions.
\item The functions $F_n$ satisfy the differential equations
\begin{equation*}
%\label{diffintro}
d F_n(z_1,\ldots,z_n)=
\sum_{i=1}^{n-1} F_i(z_1,\ldots, z_{i}) F_{n-i}(z_{i+1},\ldots, z_n)
d\log \bigg(\frac{z_{i+1}+\cdots+z_n}{z_1+\cdots +z_i} \bigg)
\end{equation*}
which implies that the functions $f_\alpha$ satisfy
\begin{equation}
\label{diffintrotwo}
d f_\alpha=
\sum_{\substack{\beta,\gamma\in K_{>0}(\A)\\ \beta+\gamma=\alpha}}
[f_\beta, f_\gamma]\, d\log {\gamma}.
\end{equation}
\end{enumerate}
\end{thm}
The remarkable point is that the specific jumping behaviour of the classes
of semistable objects leads to universal holomorphic functions satisfying an
interesting system of non--linear differential equations.

\subsection{}\label{ss:Stokes from now}
%--------------

The main new idea of this paper is that a stability condition $Z$ on $\A$  can naturally
be interpreted as defining {\it Stokes data} for an irregular connection on $\IP^1$ with
values in the Ringel--Hall Lie algebra $\wh{\bb}(\A)$. The discontinuous nature of the
classes of semistables as $Z$ varies corresponds to the  discontinuous behaviour of
the Stokes factors of an isomonodromic family of irregular connections as the Stokes
rays collide and separate. Moreover, Joyce's holomorphic functions $f_\alpha$ on
$\Stab(\A)$ can be interpreted as defining the residues of this family of connections.

\subsection{}
%--------------

Our starting point was the observation that the differential equation \eqref
{diffintrotwo} has the same form as the equation for isomonodromic deformations
of irregular connections on $\IP^1$ written down for the group $GL_n(\IC)$
by Jimbo--Miwa--Ueno \cite{JMU} and extended to an arbitrary complex,
reductive Lie group $G$ by Boalch \cite[Lemma 16]{Boalch2}.

In more detail, let $\g$ be the Lie algebra of $G$ and fix a Cartan subalgebra $\h
\subset\g$. Let $\Phi\subset \h^*$ be the associated root system and $\g=\h\oplus
\god$, with $\god=\bigoplus_{\alpha\in\Phi}\g_\alpha$ the corresponding root space
decomposition. Define a meromorphic connection on the trivial $G$--bundle over
$\IP^1$ by
\begin{equation}
\label{start}
\nabla=d-\bigg(\frac{Z}{t^2}+\frac{f}{t}\bigg)dt,
\end{equation}
where $f=\sum_{\alpha\in \Phi} f_\alpha\in\god$ and $Z$ is a regular element of
$\h$. The connection $\nabla$ has a pole of order 2 at the origin and a pole of
order 1 at infinity.

The gauge equivalence class of such a connection is determined by its Stokes data.
This data consists of a set of {\it Stokes rays}, namely the subsets of $\C$ of the form
$\R_{>0}Z(\alpha)$ for $\alpha\in\Phi$ and, for each such ray $\ell$ a corresponding
{\it Stokes factor} $S_\ell\in G$. As $Z$ varies the Stokes rays move, but if the element
$f\in\god$ evolves according to the differential equation
\begin{equation}
\label{diffintrothree}
d f_\alpha=
\sum_{\substack{\beta,\gamma\in\Phi\\ \beta+\gamma=\alpha}}
[f_\beta, f_\gamma]\, d\log {\gamma},
\end{equation}
then the Stokes factors are locally constant, and when two rays collide or separate
the product of the corresponding Stokes factors remains constant. Such deformations
of $\nabla$ are called {\it isomonodromic}.

\subsection{}
%--------------
\label{pie}

The striking similarity of the differential equations \eqref{diffintrotwo} and
\eqref{diffintrothree} suggests that that the classes $\epsilon_\alpha$
introduced by Joyce should be regarded as logarithms of Stokes factors
for a connection of the form \eqref{start}.

This interpretation is further corroborated by the following result of Reineke
\cite{Reineke}. Since $\A$ has finite length, all stability conditions $Z$ on
$\A$ have the \HN property: any non--zero module $M$ has a unique filtration
$$0=M_0\subset M_1\subset\cdots\subset M_n=M$$
where the successive factors $F_i=M_i/M_{i-1}$ are $Z$--semistable and
of strictly descending phases: $\phi(F_1)>\cdots>\phi(F_n)$.

Using the product in $\H(\A)$, this readily translates into the following identity
\cite[Proposition 4.8]{Reineke}
%can be written simply in the form
\begin{equation}\label{eq:Reineke}
\kappa_{\gamma}=
\sum_{n\geq 1}
\sum_{\substack{\gamma_1+\cdots+\gamma_n=\gamma\\
\phi(\gamma_1)>\cdots>\phi(\gamma_n)}}
\delta_{\gamma_1}*\cdots*\delta_{\gamma_n},
\end{equation}
where $\kappa_\gamma$ and $\delta_\alpha$ are the elements defined in \ref
{ss:Hall} and \ref{ss:Joyce review}.

Reineke's equation \eqref{eq:Reineke} can be rewritten in a more suggestive form
as follows. Given a ray $\ell=\R_{>0}\exp(i\pi\phi)\subset\C^*$ with $\phi\in (0,1)$, let
$\SS_\ell$ be the characteristic function of all semistable modules of phase $\phi$
(here we include the zero module).
Similarly, let $\Splus$ be the function which is equal to 1 on all modules\footnote{The
function $\Splus$ should not be confused with the identity element of the Hall algebra
$\H(\A)$ which is the characteristic function $1_0$ of the zero module.}. These functions
define elements in the completion $\CA$ of $C(\A)$ with respect to its $K_{\geq 0}(\A)
$--grading, and lie in fact in the pro--unipotent group $\wh{N}(\A)$ of invertible
grouplike elements in $\CA$ whose Lie algebra is $\wh{\n}(\A)$. The relation
\eqref{eq:Reineke} may then be rewritten as the following identity in $\wh{N}(\A)$
$$\Splus=\stackrel{\curvearrowright}{\prod_\ell}\SS_{\ell},$$
%$$S_+=\clockwise_\ell S_{\ell},$$
where the product over rays is taken in clockwise order.

The above equation is precisely what expresses the Stokes multiplier of a
connection of the form \eqref{start} relative to the upper half--plane $\IH$
in terms of its Stokes factors $S_\ell$.
The analogy with Stokes phenomena proceeds further: there is a countable
set of rays $\R_{>0} Z(\alpha)$ for $\alpha\in K_{>0}(\A)$, which can collide
or separate as $Z$ varies, and to each such ray $\ell$ is associated an element
$\SS_\ell$ of a group with the property that the ordered product of these
elements $\SS_\ell$ remains constant.

\subsection{}
%--------------

To make these analogies precise, we study in \cite{BTL} irregular
connections of the form \eqref{start} with structure group an arbitrary
complex, affine {\it algebraic} group $G$ so as to encompass the
pro--solvable group $\wh{B}(\A)$ corresponding to the \RH Lie algebra
$\wh{\bb}(\A)$. We also explicitly solve a Riemann--Hilbert problem for
these connections by expressing their residue $f$ at $t=0$ in terms of
their Stokes data.

We rely on these results to prove in this paper that the characteristic
functions $\SS_\ell$ of semistable objects of a given phase {\it are}
the Stokes factors of a unique connection of the form \eqref{start}.
The residue $f$ of this connection is given precisely by Joyce's
generating function \eqref{liedseries}. Moreover, as the stability
condition varies, the connection varies isomonodromically, thus
leading to a natural derivation of Joyce's PDE as an isomonodromic
deformation equation.

More precisely, let $\wh{P}$ be the holomorphically trivial, principal
$\wh{B}(\A)$--bundle on $\IP^1$. Let $Z\in\Stab(\A)\subset\h(\A)$ be
a stability condition and consider connections on $\wh{P}$ of the form
\begin{equation}\label{eq:the form}
%\nabla_{\A,Z}=
\nabla=d-\left(\frac{Z}{t^2}+\frac{f}{t}\right)dt
\end{equation}
where $f\in\wh{\n}(\A)$. Our main result is the following

\begin{thm}\label{th:le main intro}\hfill
\begin{enumerate}
\item
There exists a unique connection $\nabla_{\A,Z}$ of the form
\eqref{eq:the form} with Stokes data given by either of the following
equivalent conditions:
\begin{enumerate}
\item the Stokes factor corresponding to a Stokes ray $\ell=
\IR_{>0}\exp(i\pi\phi)$ is the characteristic function $\SS_\ell$
of $Z$--semistable modules of phase $\phi$.
\item The Stokes multipliers $S_+,S_-$ relative to the ray $r=
\IR_{>0}$ are the function $\Splus$ which takes the value 1
on every module, and the identity element $1_0$ respectively.
\end{enumerate}
\item The components of the residue $f=\sum_{\alpha\in K
_{>0}(\A)}f_\alpha$ of $\nabla_{\A,Z}$ at $0$ are given by
a Lie series in the elements $\{\epsilon_\beta\}$ of the form
$$f_\alpha=
\sum_{n\geq 1}
\sum_{\substack{\alpha_i\in K_{>0}(\A)\\
\alpha_1+\cdots+\alpha_n=\alpha}}
J_n(Z(\alpha_1),\ldots, Z(\alpha_n))\,
\epsilon_{\alpha_1}*\cdots*\epsilon_{\alpha_n}.$$
where $J_n:(\IC^*)^n\to\IC$ are holomorphic functions with branchcuts
which coincide with Joyce's functions $F_n$ on the domain where they
are holomorphic.

\item
As $Z$ varies in $\Stab(\A)$, the family of connections $\nabla_{\A,Z}$
varies isomonodromically. In particular, $f_\alpha(Z)$ is a holomorphic
function of $Z$ and satisfies the PDE
$$df_\alpha=
\sum_{\beta+\gamma=\alpha}[f_\beta,f_\gamma]d\log\gamma$$
\end{enumerate}
\end{thm}

\subsection{}
%--------------

We conclude with a detailed description of the contents of this paper. In Section
\ref{se:Irregular connections}, we review the definition of the Stokes data of an
irregular connection. In Section \ref{se:Stokes map}, we state the results of \cite
{BTL} on the computation of the corresponding Stokes map and the Taylor series
of its inverse in terms of multilogarithms. In Section \ref{se:Hall}, we review the
construction of the \RH algebra $\H(\A)$ of an abelian category $\A$ following
\cite{Joyce1}. In Section \ref{se:stability conditions}, we explain Joyce's construction
of $\H(\A)$--valued invariants which count semistable objects in $\A$. Section
\ref{se:stability and Stokes} contains our main result. We show that a stability
condition $A$ on $\A$ defines Stokes data for an irregular connection on $\IP^1$
with values in the \RH Lie algebra of $\A$ which varies isomonodromically with $Z$.

\subsection*{Acknowledgements} This project was begun while the first
author visited Northeastern University in March of 2007. It was completed
while the second author was a member of the Institute for Advanced Study
in Princeton during the Academic year 2007--8. He is extremely grateful
to the Institute for its financial support, hospitality
and marvelous working conditions. We would very much like to thank
Phil Boalch, Pierre Deligne, Dominic Joyce, Nikita Markarian and
Simon Willerton for useful discussions during the writing of this
paper. We are also grateful to M. Kontsevich and Y. Soibelman for
sending us a preliminary version of \cite{KS} and to D.
Huybrechts, B. Leclerc and J. Schroer for pointing out small
inaccuracies in an earlier version of this paper.

% **********************************************************************************************************
% **********************************************************************************************************
% **********************************************************************************************************
% **********************************************************************************************************
% **********************************************************************************************************

\section{Irregular connections and Stokes phenomena}
\label{se:Irregular connections}
%==========================================

We review in this section the definition of Stokes data for irregular
connections on $\IP^1$. Our exposition follows \cite{BTL} which,
in turn, is patterned on \cite{Boalch3,Boalch2}. Unlike \cite{Boalch3,
Boalch2} and earlier treatements however (see, \eg \cite{BJL}), we
do not restrict ourselves to connections whose structure group is
reductive. We consider instead the case of an arbitrary algebraic
group since this larger class encompasses the pro--solvable \RH
groups of abelian categories.

\subsection{Algebraic groups} %Algebraic groups
%-----------------------------------

By an {\it algebraic group}, we shall always mean an affine algebraic
group $G$ over $\IC$. By a {\it \fd representation} of $G$, we
shall mean a rational representation, that is a morphism $G\to
GL(V)$, where $V$ is a \fd complex vector space. An algebraic
group always possesses a faithful \fd representation and may
therefore be regarded as a {\it linear} algebraic group, that is
a (Zariski) closed subgroup of some $GL(V)$, see \eg \cite{Hu}.

Let $\g$ be the Lie algebra of $G$. If $\rho:G\to GL(V)$
is a \fd representation, we %abusively
denote its differential $\g\to \gl(V)$ by the same symbol.

\subsection{}\label{ss:setup}
%--------------

Let $G$ be an algebraic group, $H\subset G$ a maximal
torus in $G$ and $\g,\h$ their Lie algebras. The following are
the examples which will be most relevant to us
\begin{enumerate}
\item $G=GL_n(\IC)$ and $H$ is the torus consisting of diagonal,
invertible matrices.
\item $G$ is a complex, semisimple Lie group and $H\subset G$
is a maximal torus.
\item $G=H\ltimes N$, where $H$ is a torus acting on a unipotent
group $N$.
\end{enumerate}
As outlined in the Introduction, and further explained in Section \ref
{se:Hall}, case (iii) arises naturally when studying the abelian category
$\A=\Mod(R)$ of \fd representations of a \fd algebra $R$ over $\IC$.
In that case, $N=\wh{N}(\A)$ is the (pro--)unipotent group whose Lie
algebra is the Ringel--Hall Lie algebra $\nA$ of $\A$ and $H$
is the torus whose character lattice is the Grothendieck group
$K(\A)$.

\subsection{}
%--------------

Let $\bfX(H)$ be the group of characters of $H$ and $\bfX(H)\cong
\Lambda\subset\h^*$ the lattice spanned by the differentials of elements
in $\bfX(H)$. For any $\lambda\in\Lambda$ we denote the unique
element of $\bfX(H)$ with differential $\lambda$ by $e^{\lambda}$.

Decompose $\g$ as
\begin{equation}\label{eq:root decomposition 1}
\g=\h\oplus\god\qquad\text{with}\qquad
\god=\bigoplus_{\alpha\in\PhiZ}\g_\alpha,
\end{equation}
where $\Phi\subset\Lambda\setminus\{0\}$ is a finite subset and $H$
acts on $\g_\alpha$ via the character $e^\alpha$. We refer to the elements
of $\Phi$ as the {\it roots of $G$}. Since $H$ is a maximal torus, $\Phi$ is
independent of the choice of $H$.

\subsection{The irregular connection $\mathbf{\nabla}$}\label{ss:assumptions}
%--------------------------------------------------------------------

Let $P$ be the holomorphically trivial, principal $G$--bundle on
$\IP^1$ and consider the meromorphic connection on $P$ given
by
\begin{equation}
\label{nab}
\nabla=d-\bigg(\frac{Z}{t^2}+\frac{f}{t}\bigg)dt.
\end{equation}
where $Z,f\in\g$.

The connection $\nabla$ has a pole of order 2 at $t=0$ and a pole of
order 1 at $\infty$. We henceforth make the following assumptions:
\begin{enumerate}
\item[$(Z)$] $Z\in\hreg=\h\setminus\bigcup_{\alpha\in\Phi}\Ker(\alpha)$
is a regular element of $\h$. % semisimple.
\item[($f$)] $f\in\god\subset\g.$%:=\bigoplus_{\alpha\in\PhiZ}\g_\alpha$.
%The projection of $f$ onto $\h$ corresponding to the
%decomposition \eqref{eq:root decomposition 1} is zero.
\end{enumerate}
The reader unfamiliar with algebraic groups may wish to consider the
case when $G=\GL_n(\C)$ and $Z$ is a diagonal matrix with distinct
eigenvalues. Condition ($f$) is then the statement that the diagonal
entries of the matrix $f$ are zero.

\subsection{Stokes rays and sectors}
%--------------------------------------------

\begin{defn}\label{de:Stokes rays}
A {\it ray} is a subset of $\IC^*$ of the form $\R_{>0}\exp(i\pi\phi)$.
The {\it Stokes rays} of the connection $\nabla$ are the rays $\R_{>0}Z
(\alpha)$, $\alpha\in\PhiZ$. The {\it Stokes sectors} are the open regions
of $\C^*$ bounded by them. A ray is called {\it admissible} if it is not a
Stokes ray.
\end{defn}

\subsection{Canonical fundamental solutions}
%-------------------------------------------------------

The Stokes data of the connection $\nabla$ are defined using fundamental
solutions with prescribed asymptotics. We first recall how these are
characterised.

Given a ray $r$ in $\C$, we denote by $\IH_r$ the corresponding
half--plane
\begin{equation}\label{eq:halfplane}
\IH_r=\{z=uv:u\in r, \Re(v)>0\}\subset \C^*.
\end{equation}
The following basic result is well--known for $G=\GL_n(\IC)$ and $Z$
regular (see, \eg \cite[pp. 58--61]{Wasow}) and was extended in \cite
{Boalch2} to the case of complex reductive groups. It is proved in \cite
{BTL} for an arbitrary algebraic group.

\begin{thm}
\label{jurk}
Given an admissible ray $r$, there is a unique holomorphic function
$Y_r:\IH_r\to G$ such that
\begin{gather}
\frac{dY_r}{dt}=\left(\frac{Z}{t^2}+\frac{f}{t}\right)Y_r
\label{eq:diff equ}\\
Y_r\cdot e^{Z/t}\to 1\quad\text{as}\quad\text{$t\to 0$ in $\IH_r$}
\label{eq:asym}
\end{gather}
\end{thm}

\subsection{}
%--------------

The uniqueness part of Theorem \ref{jurk} relies upon the following
statement which is proved in \cite{BTL} (see \cite[Lemma 22]{Boalch2}
for the case of $G$ reductive).

\begin{prop}\label{pr:spectral}
Let $r,r'$ be two rays such that $r\neq -r'$, and $g\in G$ an element
such that
$$e^{-Z/t} \cdot g\cdot e^{Z/t}\to 1\quad\text{ as }\quad t\to 0\text{ in }\IH_r\cap \IH_{r'}.$$
Then, $g=\exp(X)$ where $X$ lies in
$$\bigoplus_{\alpha:Z(\alpha)\in\ol{\Sigma}(r,r')}\g_\alpha\subset\g,$$
with $\ol{\Sigma}(r,r')\subset\IC^*$ the closed convex sector bounded
by $r$ and $r'$.
\end{prop}

\noindent
Proposition \ref{pr:spectral} implies in particular that if the rays $r,r'$
are admissible and such that the sector $\ol{\Sigma}(r,r')$ does not
contain any Stokes rays of $\nabla$, the element $g\in G$ determined by
$$Y_r(t)=Y_{r'}(t)\cdot g\quad\text{for $t\in\IH_r\cap\IH_{r'}$}$$
is equal to 1.

\subsection{Stokes factors}\label{ss:Stokes factors}
%--------------------------------

Assume now that $\ell$ is a Stokes ray. Let $r_\pm$ be small clockwise
(resp. anticlockwise) perturbations of $\ell$ such that the convex sector
$\ol{\Sigma}(r_-,r_+)$ does not contain any Stokes rays of $\nabla$
other than $\ell$. 

\begin{defn}\label{de:Stokes factor}
The {\it Stokes factor} $S_\ell$ corresponding to $\ell$ is the element of
$G$ defined by
$$Y_{r_+}(t)= Y_{r_-}(t)\cdot S_\ell\text{ for }t\in \IH_{r_+}\cap \IH_{r_-}.$$
\end{defn}

\noindent By Proposition \ref{pr:spectral}, the definition of $S_\ell$ is
independent of the choice of $r_\pm$.

\subsection{Stokes multipliers}\label{ss:multipliers}
%------------------------------------

An alternative but closely related system of invariants are the Stokes
multipliers of the connection $\nabla$. These depend on a choice of
a ray $r$ such that both $r$ and $-r$ are admissible.

\begin{defn}
The {\it Stokes multipliers} of $\nabla$ corresponding to $r$ are the
elements $S_\pm\in G$ defined by
\[ Y_{r,\pm}(t)=Y_{-r}(t)\cdot S_\pm,\quad t\in\IH_{-r}\]
where $Y_{r,+}$ and $Y_{r,-}$ are the analytic continuations of $Y_r$
to $\IH_{-r}$ in the anticlockwise and clockwise directions respectively.
\end{defn}

\noindent
Note that the multipliers $S_\pm$ remain constant under a perturbation
of $r$ so long as $r$ and $-r$ do not cross any Stokes rays.

\subsection{}
%--------------

To relate the Stokes factors and multipliers, set $r=\R_{>0}\exp(i\pi\theta)$
and label the Stokes rays as $\ell_j=\R_{>0}\exp(i\pi\phi_j)$, with $j=1,\ldots,
m_1+m_2$ where
$$\theta<\phi_1<\cdots<\phi_{m_1}< \theta+1<\phi_{m_1+1}<\cdots<\phi_{m_1+m_2}<
\theta+2$$
The following result is immediate upon drawing a picture
\begin{lemma}
\label{veryeasy}
The following holds
$$S_+=S_{\ell_{m_1}}  \cdots S_{\ell_1}\aand
S_-=S_{\ell_{m_1+1}}^{-1}\cdots S_{\ell_{m_1+m_2}}^{-1}$$
\end{lemma}

The Stokes factors therefore determine the Stokes multipliers for
any ray $r$. In fact, conversely, the Stokes multipliers for a single
ray $r$ determine all the Stokes factors, although this is not so
easy to see. It will follow from Proposition \ref{pr:Stokes HN Reineke}
below.

\subsection{Isomonodromic families of connections}\label{ss:IM}
%---------------------------------------------------------------

We discus below isomonodromic deformations of $\nabla$. Their main
interest from our point of view, indeed one of the starting points of the
present work, is the isomonodromy equations \eqref{iso} which bear a
striking resemblance to the non--linear system of PDEs \eqref{diffintrotwo}
appearing in Joyce's work \cite{Joyce}.

Let $\U\subset\hreg$ be an open set and consider a family of connections
of the form
$$\nabla(Z)=d-\bigg(\frac{Z}{t^2}+\frac{f(Z)}{t}\bigg)dt$$
where $Z$ varies in $\U$ and the dependence of $f(Z)\in\god$ \wrt $Z$
is arbitrary.

\begin{defn}
The family of connections $\nabla(Z)$ is {\it isomonodromic} if for any
$Z_0\in\U$, there exists a neighborhood $Z_0\in\U_0\subset\U$ and
a ray $r$ such that $\pm r$ are admissible for all $\nabla(Z)$, $Z\in\U
_0$ and the Stokes multipliers $S_\pm(Z)$ of $\nabla(Z)$ relative to
$r$ are constant on $\U_0$.
\end{defn}

The isomonodromy of the family $\nabla(Z)$ may also be defined as the
constancy of the Stokes factors. This requires a little more care since,
as pointed out in \cite[pg. 190]{Boalch} for example, Stokes rays may
split into distinct rays under arbitrarily small deformations of $Z$.

Call a sector $\Sigma\subset\IC^*$ {\it admissible} if its boundary rays
are admissible.

\begin{prop}
The family of connections $\nabla(Z)$ is isomonodromic if, and only
if, for any connected open subset $\U_0\subset \U$ and any convex
sector $\Sigma$ which is admissible for all $\nabla(Z)$, $Z\in \U_0$,
the clockwise product
$$\clockwise_{\ell\subset\Sigma}S_\ell(Z)$$
of Stokes factors corresponding to the Stokes rays contained in
$\Sigma$ is constant on $\U_0$.
\end{prop}
\begin{pf}
This follows from the fact that Stokes factors and multipliers determine
each other by Lemma \ref{veryeasy} and Proposition \ref{pr:Stokes HN Reineke}.
\end{pf}

\subsection{Isomonodromy equations}
%----------------------------------------------

The following characterisation of isomonodromic deformations was obtained
by Jimbo--Miwa--Ueno in the case $G=GL_n$ \cite{JMU} and adapted to the
case of a complex, reductive group by Boalch \cite[Appendix]{Boalch}.
Its proof carries over verbatim to the case of a general algebraic group.

\begin{thm}\label{th:IMD}
Assume that $f$ varies holomorphically in $Z$. Then, family of connections
$\nabla(Z)$ is isomonodromic if, and only if $f$ satisfies the PDE
\begin{equation}
\label{iso}
d f_\alpha=\sum_{\substack{\beta,\gamma\in\Phi\\\beta+\gamma=\alpha}}
[f_\beta, f_\gamma]\, d\log\gamma.
\end{equation}
\end{thm}

\remark\label{rk:IMD} The equations \eqref{iso} form a first order system
of {\it integrable} non--linear PDEs and therefore possess a unique
holomorphic solution $f(Z)$ defined in a neighboorhood of a fixed
$Z_0\in\hreg$ and subject to the initial condition $f(Z_0)=f_0\in[\h,\g]$.

\remark Jimbo--Miwa--Ueno and Boalch also give an alternative characterisation
of isomonodromy in terms of the existence of a flat connection on $\IP^1\times U$
which has a logarithmic singularity on the divisor $\{t=\infty\}$ and a pole of order
2 on $\{t=0\}$, and restricts to $\nabla(Z)$ on each fibre $\{Z\}\times\IP^1$. This
connection is given by
\[\ol{\nabla} = d - \bigg[\bigg(\frac{Z}{t^2}+\frac{f}{t}\bigg)dt
 +  \sum_{\alpha\in \Phi} f_\alpha\frac{d\alpha}{\alpha}
+\frac{ dZ}{t}\bigg].\]
One can check directly that the flatness of this connection
is equivalent to \eqref{iso}.

\section{The Stokes map}\label{se:Stokes map}
%===================

In this section we state the main results of \cite{BTL}. These
express the logarithms of the Stokes factors of the connection
$\nabla$ as explicit, universal Lie series in the variables $f_\alpha$
and, conversely, the $f_\alpha$ as Lie series in the logarithms
of the Stokes factors, thus explicitly solving a Riemann--Hilbert
problem.

\subsection{Completion with respect to \fd representations}
%------------------------------------------------------------------------

Our formulae are more conveniently expressed inside the completion
$\wh{U\g}$ of $U\g$ \wrt the \fd representations of $G$. We review
below the definition of $\wh{U\g}$.

Let $\Vec$ be the category of \fd complex vector spaces and $\RRep
(G)$ that of \fd representations of $G$. Consider the forgetful functor
$$F:\RRep(G)\rightarrow\Vec.$$
By definition, $\wh{U\g}$ is the algebra of endomorphisms of $F$.
Concretely, an element of $\wh{U\g}$ is a collection $\Theta=\{\Theta
_V\}$, with $\Theta_V\in\End_{\IC}(V)$ for any $V\in\RRep(G)$, such
that for any $U,V\in\RRep(G)$ and $T\in\Hom_G(U,V)$, the following
holds
$$\Theta_V\circ T=T\circ\Theta_U$$

There are natural homomorphisms $U\g\to\wh{U\g}$ and $G\to\wh{U\g}$
mapping $x\in U\g$ and $g\in G$ to the elements $\Theta(x)$, $\Theta(g)$
which act on a \fd representation $\rho:G\to GL(V)$ as $\rho(x)$ and $\rho
(g)$ respectively. These homomorphisms are well--known to be injective
(see \eg \cite[Lemma 4.1]{BTL}) and we will use them to think of $U\g$
as a subalgebra of $\wh{U\g}$ and $G$ as a subgroup of the group of
invertible elements of $\wh{U\g}$ respectively.

\subsection{Representing Stokes factors}
%-------------------------------------------------
\label{rep}

Fix a Stokes ray $\ell$ of the connection $\nabla$. We show below how to
represent the corresponding Stokes factor $S_\ell$ in two different ways:
by elements $\epsilon_\alpha\in\god$ and by elements $\delta_\gamma\in
U\g$.

Consider the subalgebra
$$\n_\ell=\bigoplus_{\alpha:Z(\alpha)\in\ell} \g_\alpha\subset \g.$$
The elements of $\n_\ell$ are nilpotent, that is they act by nilpotent
endomorphisms on any \fd representation of $G$. It follows that the
exponential map $\exp\colon\n_\ell\to G$ is an isomorphism onto the
unipotent subgroup $N_\ell=\exp(\n_\ell)\subset G$.

By Proposition \ref{pr:spectral}, the Stokes factor $S_\ell$ lies in
$N_\ell$. For the first representation of $S_\ell$, write
\begin{equation}
\label{bubbly}
S_\ell=\exp\bigg(\sum_{\alpha:Z(\alpha)\in\ell}\epsilon_\alpha\bigg)
\end{equation}
for uniquely defined elements $\epsilon_\alpha\in\g_\alpha$.

For the second, we compute the exponential
\eqref{bubbly} in $\wh{U\g}$ and decompose the result along the
weight spaces
$$\wh{U\g}_\gamma=
\{x\in\wh{U\g}|\ad(h)x=\gamma(h)x,\;\forall h\in\h\},\;\gamma\in\h^*$$
of the adjoint action of $\h$. This yields elements $\delta_\gamma
\in (U\n_\ell)_\gamma$ such that
\begin{equation}
\label{hold}
S_{\ell}
=1+\sum_{\gamma\in\LambdaZ:Z(\gamma)\in\ell}\delta_\gamma,
\end{equation}
where $\LambdaZ\subset\h^*$ is the lattice generated by the set of roots
$\PhiZ$ and the above identity is to be understood as holding in any \fd
\rep of $G$.

These two representations of $S_\ell$ are related as follows.

\begin{lemma}\label{le:epsilon delta}\hfill
\begin{enumerate}
\item
Let $\gamma\in\LambdaZ$ be such that $Z(\gamma)$ lies on the Stokes
ray $\ell$. Then, $\delta_\gamma$ is given by the finite sum
\begin{equation}
\label{exp}
\delta_\gamma=
\sum_{n\geq 1}
\sum_{\substack{\alpha_i\in\PhiZ\\Z(\alpha_i)\in\ell,\\\alpha_1+\cdots+\alpha_n=\gamma}}
\frac{1}{n!}\;\epsilon_{\alpha_1} \cdots \epsilon_{\alpha_n}.
\end{equation}
\item
Conversely, let $\alpha\in\PhiZ$ be such that $Z(\alpha)\in\ell$. Then,
$\epsilon_\alpha$ is given by the finite sum
\begin{equation}
\label{log}
\epsilon_\alpha=
\sum_{n\geq 1}
\sum_{\substack{\gamma_i\in\LambdaZ\\Z(\gamma_i)\in\ell,
\\\gamma_1+\cdots+\gamma_n=\alpha}}
\frac{(-1)^{n-1}}{n}\,\delta_{\gamma_1}\cdots\delta_{\gamma_n}.
\end{equation}
\end{enumerate}
\end{lemma}
\begin{pf}
These are the standard expansions of $\exp\colon\n_\ell\to N_\ell$ and
$\log\colon N_\ell\to\n_\ell$.
\end{pf}

%Once again we can assemble the elements $\delta_\gamma$ corresponding to
%different Stokes rays and form the sum\comment{Any reason we might want to
%make $\delta_0$ equal to 1 so that, for example, $S=\delta$? The question
%applies to Joyce's stuff as well.}
%$$\delta=\sum_{0\neq \gamma\in\LambdaZ} \delta_\gamma\in \wh{U\g}.$$

\subsection{Formula for the Stokes factors}
%----------------------------------------------------

We now give an explicit formula for the Stokes factors of the connection
$\nabla$ in terms of iterated integrals.

\begin{defn}\label{de:M}
Set $M_1(z_1)=2\pi i$ and, for $n\geq 2$, define the function $M_n:
(\IC^*)^n\to\IC$ by the iterated integral
$$M_n(z_1,\ldots,z_n)=
2\pi i \big.\int_{C}\frac{dt}{t-s_1}\circ \cdots \circ \frac{dt}{t-s_{n-1}},$$
where $s_i=z_1+ \cdots +z_i$, $1\leq i\leq n$ and the path of integration $C$
is the line segment $[0,s_n]$, perturbed if necessary to avoid any point $s_i\in
[0,s_n]$ by small clockwise arcs.
\end{defn}

\begin{thm}[\cite{BTL}]\label{one2}
The weight components $\delta_\gamma$ of the Stokes factor $S_\ell$
corresponding to the ray $\ell$ are given by
\begin{equation}
\label{mainy}
\delta_\gamma=
\sum_{n\geq 1}
\sum_{\stackrel{\alpha_i\in\PhiZ}{\alpha_1+\cdots+\alpha_n=\gamma}}
M_n(Z(\alpha_1),\ldots, Z(\alpha_n))\,
f_{\alpha_1} f_{\alpha_2}\cdots f_{\alpha_n},
\end{equation}
where the equality is to be understood as holding in any \fd \rep of $G$
and the sum over $n$ is absolutely convergent.
\end{thm}

\subsection{The Stokes map}\label{ss:Stokes}
%-----------------------------------

Since the sets $\{\alpha\in\Phi:Z(\alpha)\in\ell\}$ partition $\PhiZ$ as $\ell$
ranges over the Stokes rays of $\nabla$, we may assemble the elements
$\epsilon_\alpha$ corresponding to different Stokes rays and form the sum
$$\epsilon=
\sum_{\alpha\in \PhiZ}\epsilon_\alpha\in\bigoplus_{\alpha\in\PhiZ}\g_\alpha.$$
For fixed $Z\in\h$, we shall refer to the map
$$\calS:\bigoplus_{\alpha\in\PhiZ}\g_\alpha\longrightarrow
\bigoplus_{\alpha\in\PhiZ}\g_\alpha$$
mapping $f$ to $\epsilon$ as the {\it Stokes map}.

\subsection{Formula for the Stokes map}
%--------------------------------------------------

We next state a formula for the Stokes map giving the element $\epsilon$ in
terms of $f$. We first define  the special functions appearing in this formula.

\begin{defn}
The function $L_n:(\IC^*)^n\to\IC$ is given by $L_1(z_1)=2\pi i$ and, for
$n\geq 2$,
$$L_n(z_1,\ldots,z_n)=
\sum_{k=1}^n\;\sum_{\substack{
0=i_0<\cdots<i_k= n\\
s_{i_j}-s_{i_{j-1}}\in\IR_{>0}\cdot s_n}}
\frac{(-1)}{k}^{k-1}\;
\prod_{j=0}^{k-1} M_{i_{j+1}-i_j}(z_{i_j +1} ,\ldots,z_{i_{j+1}}),$$
where $s_j=z_1+\cdots+z_j$.
\end{defn}

\begin{remark}
Note that on the open subset
\[(z_1,\ldots,z_n) \in (\C^*)^n\text{ such that } s_i\notin [0,s_n]\text{ for }0<i<n\]
the inner sum above is empty unless $k=1$ and one therefore has
$$L_n(z_1,\ldots,z_n)=M_n(z_1,\ldots,z_n).$$
Thus $L_n$ agrees with $M_n$ on the open subset where it is holomorphic and
differs by how it has been extended onto the cutlines.
\end{remark}

\subsection{}
%--------------

The functions $L_n$ are more complicated to define than the functions $M_n$.
Unlike the latter however, they give rise to Lie series as the following result shows.

\begin{thm}[\cite{BTL}]\label{one}\hfill
\begin{enumerate}
\item Let $x_1,\ldots, x_n$ be elements in a Lie algebra $\L$. For
any $(z_1,\ldots,z_n)\in (\C^*)^n$, the finite sum
$$\sum_{\sigma\in \operatorname{Sym}_n}
L_n(z_{\sigma(1)}, \cdots, z_{\sigma(n)})
x_{\sigma(1)} \cdots x_{\sigma(n)}$$
is a Lie polynomial in $x_1,\ldots,x_n$ and therefore lies in $\L\subset U\L$.
\item The element $\epsilon=\calS(f)$ is given by a Lie series in the variables
$\{f_\alpha\}_{\alpha\in\PhiZ}$, given by
\begin{equation}
\label{main}
\epsilon_\alpha=
\sum_{n\geq 1}\sum_{\stackrel{\alpha_i\in\PhiZ}{\alpha_1+\cdots+\alpha_n=\alpha}}
L_n(Z(\alpha_1),\ldots, Z(\alpha_n)) f_{\alpha_1} f_{\alpha_2}\cdots f_{\alpha_n}
\end{equation}
As a series in $n$, \eqref{main} converges uniformly on compact subsets of $\god$.
%where the sum over $n$ is absolutely convergent.
\end{enumerate}
\end{thm}

\subsection{Inverse of  the Stokes map}
%-------------- (inverse IRH map)

By Theorem \ref{one}, the Stokes map $\calS:\god\to\god$ is holomorphic,
satisfies $\calS(0)=0$ and its differential at $f=0$ is invertible. By the
inverse function Theorem, $\calS$ possesses an analytic inverse
$\calS^{-1}$ defined on a neighborhood of $\epsilon=0$.

\begin{thm}[\cite{BTL}]\label{th:Stokes inverse}
The Taylor series of $\calS^{-1}$ at $\epsilon=0$ is given by a Lie series
in the variables $\{\epsilon_\alpha\}_{\alpha\in\Phi}$ of the form
\begin{equation}
\label{main2}
f_\alpha=
\sum_{n\geq 1}
\sum_{\stackrel{\alpha_i\in\PhiZ}{\alpha_1+\cdots+\alpha_n=\alpha}}
J_n(Z(\alpha_1),\ldots, Z(\alpha_n))\,
\epsilon_{\alpha_1} \epsilon_{\alpha_2}\cdots \epsilon_{\alpha_n}
\end{equation}
for some special functions $J_n\colon (\C^*)^n \to \C$.
\end{thm}

\begin{remark}
It follows from Theorem \ref{th:Stokes inverse} that whenever the sum \eqref
{main2} is absolutely convergent over $n$, it inverts the Stokes map $\calS$,
in that the connection \eqref{nab} determined by $f=\sum_{\alpha}f_\alpha$
has Stokes factors given by \eqref{bubbly}.
\end{remark}

\subsection{The functions $\mathbf{J_n}$}
%---------------------------------------------------

The functions $J_n$ appearing in Theorem \ref{th:Stokes inverse} are
explicitly described in \cite{BTL} as sums of products of the functions
$L_n$ indexed by plane rooted trees. For example
\begin{multline*}
(2\pi i)^3 J_3(z_1,z_2,z_3)= \\ L_2(z_1,z_2)L_2(z_1+z_2,z_3)
 - L_3(z_1,z_2,z_3) + L_2(z_1,z_2+z_3) L_2(z_2,z_3)
\end{multline*}
corresponding to the three distinct plane rooted trees with 3 leaves.

\begin{thm}[\cite{BTL}]\label{js}
The function $J_n\colon (\C^*)^n\to \C$ is continuous and holomorphic
on the complement of the hyperplanes
$$H_{ij}=\{z_i+\cdots+z_j=0\},\quad 1\leq i<j\leq n$$
in the domain
$$\D_n=\{(z_1,\ldots,z_n)\in (\C^*)^n|\medspace
z_i/z_{i+1} \notin \R_{>0}
\text{ for } 1\leq i< n\}$$
Moreover, it satisfies the differential equation
$$d J_n(z_1,\ldots,z_n)=
\sum_{i=1}^{n-1} J_{i}(z_1,\ldots, z_{i}) J_{n-i}(z_{i+1},\ldots, z_n)
d\log \bigg(\frac{z_{i+1}+\cdots+z_n}{z_1+\cdots +z_i} \bigg)$$
together with the conditions $J_1(z)=1/2\pi i$ and
$$J_n(z_1,\ldots,z_n)=0\quad\text{if}\quad z_1+\cdots+z_n=0$$
for $n\geq 2$.
\end{thm}

\begin{remark}
It follows from Theorem \ref{js} that the functions $J_n$ are the same as the
functions $F_n$ appearing in Joyce's paper \cite{Joyce} and alluded to in the
Introduction, at least on the dense open subset where they are holomorphic.
\end{remark}

\subsection{Representing Stokes multipliers}\label{ss:rep kappa}
%------------------------------------------------------

We next show how to represent the Stokes multipliers $S_\pm$ by an
element $\kappa\in\wh{U\g}$.

Let $r=\R_{>0}e^{i\pi\theta}$ be the ray \wrt which $S_\pm$ are defined
and $\pm i\IH_r$ the connected components of $\IC\setminus\IR\,e^
{i\pi\theta}$. These determine a partition of $\PhiZ=\PhiZ_+\sqcup\PhiZ_-$
given by
$$\PhiZ_\pm=\{\alpha\in\PhiZ:Z(\alpha)\in\pm i\IH_r\}$$
Let $\LambdaZ_\pm\subset\h^*\setminus\{0\}$ be the cones spanned
by the linear combinations of elements in $\Phi_\pm^Z$ with coefficients
in $\IN_{>0}$. Similarly to \S \ref{rep}, it follows from Proposition \ref
{pr:spectral} that there is a unique element
$$\kappa=
\sum_{\gamma\in\LambdaZ_+\sqcup\LambdaZ_-}\kappa_\gamma\in \wh{U\g}$$
such that the Stokes multipliers $S_\pm$ are equal to
$$S_+=
1+\sum_{\gamma\in\LambdaZ_+}\kappa_\gamma, \qquad
(S_-)^{-1}=1+\sum_{\gamma\in\LambdaZ_-}\kappa_\gamma,$$

Given $\gamma\in\LambdaZ_+$, set
$$\phi(\gamma)=\frac{1}{\pi}\arg Z(\gamma)\in (\theta,\theta+1).$$
The following result gives the relation between the elements $\kappa$ and $\delta$.
\begin{prop}\label{pr:Stokes HN Reineke}
\label{prop}\hfill
\begin{enumerate}
\item For all $\gamma\in\LambdaZ_+$, there is a finite sum
\begin{equation}
\label{help}
\kappa_{\gamma}=
\sum_{n\geq 1}
\sum_{\stackrel{\gamma_1+\cdots+\gamma_n=\gamma}{\phi(\gamma_1)>\cdots>\phi(\gamma_n)}} \delta_{\gamma_1}\cdots\delta_{\gamma_n},
\end{equation}
where the sum is over elements $\gamma_i\in \LambdaZ_+$.
\item Conversely, for $\gamma\in\LambdaZ_+$
\begin{equation}
\label{reinekeinverted}
\delta_{\gamma}=
\sum_{n\geq 1}\sum_{\stackrel
{\gamma_1+\cdots+\gamma_n=\gamma}{\phi(\gamma_1+\cdots+\gamma_i)> \phi(\gamma)}}
(-1)^{n-1}\kappa_{\gamma_1}\cdots \kappa_{\gamma_n},
\end{equation}
\end{enumerate}
\end{prop}
\begin{pf}
(i) follows from substituting \eqref{hold} into the formula of Lemma \ref{veryeasy}.
(ii) follows from Reineke's inversion of formula \eqref{help} \cite[Section 5]{Reineke}.
\end{pf}

\subsection{}
%--------------
\label{dia}

The following diagram summarizes the relationships between the elements $\delta,
\epsilon$ representing the Stokes factors, the element $\kappa$ representing the
Stokes multipliers, and the element $f\in\god$ defining $\nabla$.

%\[\xymatrix{&(\epsilon) \ar@/_1pc/[ldd]_{\stackrel{\text{IRH}^{-1}}{\eqref{main2}}} \ar@/_1pc/[rr]_{\exp \eqref{exp}} %&& (\delta) \ar@/_1pc/[ll]_{\log\ \eqref{log}} \ar@/_1pc/[rr]_{\eqref{help}} %&&(\kappa)\ar@/_1pc/[ll]_{\eqref{reinekeinverted}} \\ \\
%(f) \ar@/_1pc/[uur]_{\stackrel{\text{IRH}}{\eqref{main}}} \ar@/_2pc/[uurrr]_{\eqref{mainy}}}\]

%\[\xymatrix{(f) \ar@/_4pc/[rrrr]^{\eqref{mainy}} \ar@/_1pc/[rr]_{\stackrel{\text{Stokes}}{\eqref{main}}}&&(\epsilon) \ar@/_1pc/[ll]_{\stackrel{(\text{Stokes})^{-1}}{\eqref{main2}}} \ar@/_1pc/[rr]_{\stackrel{\exp}{\eqref{exp}}} && (\delta) \ar@/_1pc/[ll]_{\stackrel{\log}{\eqref{log}}} \ar@/_1pc/[rr]_{\eqref{help}} &&(\kappa)\ar@/_1pc/[ll]_{\eqref{reinekeinverted}}
 % }\]

\[\xymatrix{(f)  \ar@/_1pc/[rr]_{\stackrel{\text{\scriptsize Stokes}}{\eqref{main}}}&&(\epsilon) \ar@/_1pc/[ll]_{\stackrel{\text{\scriptsize Stokes}^{-1}}{\eqref{main2}}} \ar@/_1pc/[rr]_{\stackrel{\text{\scriptsize exp}}{\eqref{exp}}} && (\delta) \ar@/_1pc/[ll]_{\stackrel{\text{\scriptsize log}}{\eqref{log}}} \ar@/_1pc/[rr]_{\stackrel
{\text{\scriptsize clockwise multiplication}}{\eqref{help}}} &&(\kappa)\ar@/_1pc/[ll]_{\stackrel{\text{\scriptsize Reineke inversion}}{\eqref{reinekeinverted}}}
  }\]

\medskip
\noindent These related systems of invariants will appear again in Section
\ref{se:stability conditions} in the context of stability conditions on abelian
categories.

% *****************************************************************************************************
% *****************************************************************************************************
% *****************************************************************************************************
% *****************************************************************************************************
% *****************************************************************************************************

\section{Ringel--Hall algebras}\label{se:Hall}
%======================

In this section, we review the definition of the Hall algebra of an abelian
category $\A$. We shall in fact restrict ourselves to the case where $\A
=\Mod(R)$ is the category of finite--dimensional, left modules over a fixed
finite--dimensional, associative $\C$--algebra $R$. 
It should be possible to generalize our main results to include
other abelian categories, for example categories of coherent
sheaves, but rather than working in maximal generality we prefer
to focus on a case where the underlying ideas are not clouded by
technical issues.

\subsection{The Grothendieck group $K(\A)$}
%-------------------------------------------------------

Since the category $\A=\Mod(R)$ has finite length and finitely many simple
modules $S_1,\ldots, S_N$, up to isomorphism, the Grothendieck group $
K(\A)$ is a free abelian group of finite rank generated by the  classes $[S_i]$
\[K(\A)=
\Z[S_1]\oplus \cdots \oplus \Z[S_N].\]
The positive and non--negative cones $K_{>0}(\A)\subset K_{\geq 0}(\A)\subset
K(\A)$ are defined by
$$K_{> 0}(\A)=\{[M] :  0\neq M\in \A\}\aand K_{\geq 0}(\A)=K_{> 0}(\A)\sqcup \{0\}.$$

\subsection{The Ringel--Hall algebra}
%---------------------------------------------

There are many variants of the Hall algebra of $\A$, see for example \cite{Schiff}
for a survey of Hall algebras over finite fields. We shall work over $\IC$ using
constructible functions, an idea originally due to Schofield \cite{Schofield} and
later taken up by Lusztig \cite{Lusztig} and Riedtmann \cite{Riedtmann}. The
precise construction we use was sketched by Kapranov and Vasserot \cite{KV}
and described in detail by Joyce \cite{Joyce1}.
\smallskip

Recall that a complex--valued function $f\colon X\to \C$ on a variety $X$ is
{\it constructible} if it is of the form
$$f=\sum_{i=1}^k a_i 1_{Y_i}$$
for complex numbers $a_1,\ldots, a_k$ and locally--closed subvarieties $Y_i
\subset X$. Such a function can be integrated by using the Euler characteristic
as a measure \cite{MacPherson}. By definition
\[\int_X f\dchi = \sum_{i=1}^k a_i \chi(Y_i),\]
where $\chi(Z)$ is the topological Euler characteristic of a complex variety
$Z$ endowed with the analytic topology.
%\begin{itemize}
%\item if $U\subset X$ is an open subvariety then $\chi(X)=\chi(U)+\chi(X\setminus U)$ \\
%\item $\chi(X\times Y)=\chi(X)\chi(Y)$ \\
%\item $\chi(\mathbb{A}^n)=1$ for all $n\geq 0$.\end{itemize}

Given an integer $d\geq 0$, there is an affine variety $\Rep_d$ parametrising
$R$--module structures on the vector space $\C^d$. The moduli stack $\M_d$
of $R$--modules of dimension $d$ is the quotient
$$\M_d=\Rep_d/ GL_d(\IC).$$
By definition, a constructible function on $\M_d$ is a $GL_d(\IC)$--equivariant
constructible function on the affine variety $\Rep_d$%; thus it is a finite weighted
%sum of characteristic functions of $GL_d(\IC)$--invariant locally-closed\comment
%{Technically constructible, but perhaps one can reduce to locally--closed?} subsets
%of $\Rep_d$.

We define $\H_d(\A)$ to be the space of constructible functions on $\M_d$ and
set
\begin{equation}
\label{star}
\H(\A)=\bigoplus_{d\geq 0}\H_d(\A).
\end{equation}
Note that elements of $\H(\A)$ can be thought of as functions on modules that
are constant on isomorphism classes. Given $f\in\H(\A)$, we denote its value on
a module $M$ by $f(M)$. We say that an element of $\H(\A)$ is {\it supported}
on a certain class of modules to mean that its value on all other modules is
zero.

\begin{thm}[Kapranov-Vasserot \cite{KV}, Joyce \cite{Joyce1}]
There is an associative product
\[*\colon \H(\A)\tensor \H(\A) \to \H(\A)\]
for which
\[(f_1*\cdots *f_n)(M)=\int f_1(M_1/M_0) \cdots f_n(M_n/M_{n-1})\dchi,\]
where the integral is over the variety $\Flag^n(M)$ parameterising flags
\[0=M_0\subset M_1\subset \cdots \subset M_n=M\]
of submodules of $M$ of length $n$.
The characteristic function of the zero module $1_0$ is a unit for this
multiplication.
\end{thm}
\begin{pf} The case $n=2$ is proved in \cite[Theorem 4.3]{Joyce1}. Joyce
considers the morphism $\sigma(\{1,2\})$ from the stack of short exact
sequences in $\A$ to the stack of objects of $\A$ which on geometric points
takes a short exact sequence of modules
$$0\lra A\lra M\lra B\lra 0$$
to the module $M$. Note  that this morphism induces injections on stabilizer
groups since an isomorphism of short exact sequences is determined by its
action on the middle term. Thus Joyce's pushforward map on constructible
functions is just given by integrating along the fibres of the morphism $\sigma
(\{1,2\})$ which are precisely the varieties $\Flag^2(M)$ of the statement.
The extension to the product of $n$ elements follows by an induction argument.
\end{pf}

Note that the algebra $\H(\A)$ is graded by $K_{\geq 0}(\A)$:
$$\H(\A)=\bigoplus_{\gamma\in K_{\geq 0}(\A)}\H_\gamma(\A)$$
where $\H_\gamma(\A)$ is the subspace of functions supported on modules
of class $\gamma$. This grading is a refinement of the $\Z_{\geq 0}$--grading
given by \eqref{star} via the homomorphism $K_{\geq 0}(\A)\to\Z_{\geq 0}$
mapping $\gamma$ to the dimension of the modules of class $\gamma$.

\subsection{The bialgebra $C(\A)$}
%-----------------------------------------

If the moduli stack $\M_d$ has positive dimension, the space $\H_d(\A)$ of
constructible functions on it is very large since it contains the characteristic
functions of points. It is therefore usual to consider a subalgebra generated
by some natural set of elements.

For each $\gamma\in K_{\geq 0}(\A)$, let
$\kappa_\gamma\in\H_\gamma(\A)$ be the characteristic function of the
set of modules of class $\gamma$
$$\kappa_{\gamma}(M)
 =\begin{cases}
 \,1 &\text{if  $[M]=\gamma$,} \\
 \, 0 &\text{otherwise.}
\end{cases}$$
This is a constructible function because the class of a module in
$K(\A)$ is constant in families. Let $C(\A)\subset \H(\A)$ be the
subalgebra generated by the elements $\kappa_\gamma$:
$$C(\A)=
\big\langle \kappa_{\gamma}:\gamma\in K_{\geq
0}(\A)\big\rangle\subset \H(\A)$$ and note that $C(\A)$ is a
$\Z_{\geq 0}$--graded (and, {\it a fortiori}, a $K_{\geq 0}(\A)
$--graded) algebra with \fd homogeneous components.

The algebra $C(\A)$ possesses the structure of a bialgebra. To see this,
note first that the tensor product $\H(\A)\otimes\H(\A)$ embeds into $\H
(\A\times\A)$ by setting
$$(f\tensor g)(M,N)=f(M)g(N).$$
Define a map $\Delta:\H(\A)\to\H(\A\times\A)$ by
$$\Delta(f)(M,N)=f(M\oplus N)$$
The image of $\Delta$ need not be contained in $\H(\A)\otimes\H(\A)\subset\H
(\A\times\A)$ in general. The following result is due to Joyce \cite{Joyce1} and
builds upon earlier work of Ringel \cite{Ringel2}.

\begin{thm}
The map $\Delta$ restricts to a coassociative coproduct
\[\Delta\colon C(\A) \to C(\A)\tensor C(\A),\]
preserving the $K_{\geq 0}(\A)$ grading. The homomorphism $\eta\colon C(\A)\to\C$
given by evaluation on the zero module $\eta(f)=f(0)$ is a counit for $\Delta$.
The data $(*,1,\Delta,\eta)$ endows $C(\A)$ with the structure of a cocommutative
bialgebra.
\end{thm}
\begin{pf}
This follows from the proof of \cite[Theorem 4.20]{Joyce1} using the fact that
\[\Delta(\kappa_\gamma)=
\sum_{\gamma_1+\gamma_2=\gamma}
\kappa_{\gamma_1}\otimes\kappa_{\gamma_2},\]
which is immediately verified by evaluating on a pair of modules $(M,N)$.
\end{pf}

\subsection{The Ringel--Hall Lie algebra}
%-------------------------------------------------

Recall that an element $f$ of a bialgebra is {\it primitive} if $\Delta(f)=f\tensor 1+1
\tensor f$, and that the subspace of such elements is a Lie  algebra under the
commutator bracket. Recall also that a module $M\in \A$ is {\it indecomposable}
if
$$M=N\oplus P \implies N=0\text{ or }P=0.$$
In particular, the zero module is indecomposable.

\begin{lemma}\label{le:primitive}
An element $f\in C(\A)$ is primitive if, and only if it is supported on nonzero
indecomposable modules.
\end{lemma}
\begin{pf}
According to the definition of the coproduct the primitive elements of $C(\A)$
are those satisfying $f(M\oplus N)=f(M)1_{0}(N) + 1_{0}(M)f(N)$. In particular
if $M$ and $N$ are nonzero then $f(M\oplus N)=0$. Moreover $f(0)=f(0)+f(0)$
so $f(0)=0$. Hence $f$ is supported on indecomposable modules. The converse
is easily checked.
\end{pf}

We write $\n(\A)$ for the subspace of $C(\A)$ consisting of primitive elements.
Thus $\n(\A)$ is a Lie algebra which we call the {\it Ringel--Hall Lie algebra} of
$\A$. Note that the grading on $C(\A)$ induces a grading
$$\n(\A)=\bigoplus_{\alpha\in K_{>0}(\A)} \n_\alpha(\A).$$
One can use the grading  of $\n(\A)$ to form a larger Lie algebra $\bb(\A)=\h(\A)
\oplus\n(\A)$ by endowing
$$\h(\A)=\Hom_{\Z}(K(\A),\C)$$
with the trivial bracket, and setting
\[ [Z,f]=Z(\alpha) f\quad\text{for any}\quad Z\in \h(\A),f\in\n_\alpha(\A).\]
We shall  refer to $\bb(\A)$ as the {\it extended Ringel--Hall Lie algebra} of $\A$.

\begin{example}
Let $Q$ be a finite quiver and $R$ its path algebra. Assume that $Q$
does not have oriented cycles, so that $R$ is finite--dimensional. A
simple argument due to Reineke \cite[Lemma 4.4]{Reineke} shows
that $C(\A)$ coincides with the {\it composition algebra} of $\A$, that
is the subalgebra of $\H(\A)$ generated by the characteristic functions
$\kappa_{[S_i]}$ of the simple modules. In this case, $\n(\A)$ is isomorphic
to the positive part $\n_+$ of the Kac--Moody Lie algebra $\g=\n_-\oplus
\h\oplus\n_+$ corresponding to the undirected graph underlying $Q$
and $\bb(\A)$ to the corresponding Borel subalgebra $\h\oplus \n_+$.
This result was first proved over a finite field by Ringel \cite{Ringel}. A
characteristic zero result was later obtained by Schofield \cite{Schofield}.
For the exact statement made above we refer to Joyce \cite[Example 4.25]
{Joyce1} and in the finite--type case to Riedtmann \cite{Riedtmann}.
\end{example}

\subsection{Primitive generation of $C(\A)$} % Milnor-Moore
%-----------------------------------------------------

Recall that a non--zero element $c$ in a coalgebra $C$ is {\it
grouplike} if $\Delta(c)=c\otimes c$.

\begin{lemma}\label{le:grouplike}
The element $1=1_{0}$ is the only grouplike element in $C(\A)$.
\end{lemma}
\begin{pf}
If $f\in C(\A)$ is grouplike, then, for any module $M\in\A$ and $p\in\IN^*$,
$f(M^{\oplus p})=f(M)^p$. Since $f$ lies in $\H(\A)=\bigoplus_{d\in\IN}\H_
d(\A)$, it is supported on modules of dimension $\leq D$ for large enough
$D$. Choosing $p$ such that $p\dim M>D$ shows that $f(M)=0$ unless
$M=0$. In the latter case we have $f(0)=f(0\oplus 0)=f(0)^2$, whence
$f=1$ since $f\neq 0$ and therefore $f=1_{0}$.
\end{pf}

\begin{prop}
The inclusion $\n(\A)\subset C(\A)$ identifies $C(\A)$ as a bialgebra
with the universal enveloping algebra $U\n(\A)$ of $\n(\A)$.
\end{prop}
\begin{pf}
We claim first that $C(\A)$ is {\it connected}, that is that its
coradical is one--dimensional. Indeed, since $C(\A)$ is
cocommutative and defined over an algebraically closed field, any
simple subcoalgebra $C'\subset C(\A)$ is one--dimensional
\cite[page 8]{Kap} and therefore spanned by an element $c'$ which,
up to a scalar, is necessarily grouplike. By Lemma
\ref{le:grouplike}, $C'=\IC 1_{\M_0}$. The proposition now follows
from the Milnor--Moore Theorem (see, \eg \cite[thm. 21]{Kap} or
\cite[thm. 2.5.3]{Abe}).
\end{pf}

\subsection{Completion of $C(\A)$}
%------------------------------------------

For each $d\geq 1$, the subspace $C_{>d}(\A)\subset C(\A)$ of
functions supported on modules of dimension $>d$ is an ideal.
Consider the \fd algebra
$$\CbarA{d}=C(\A)/C_{>d}(\A)$$
and the corresponding inverse system
$\cdots \to \CbarA{d}\to\cdots\to\CbarA{0}=\C.$
%$$\cdots \to \CbarA{d}\to\cdots\to\CbarA{1}\to\CbarA{0}=\C.$$
By definition, the completion $\CA$ is the limit
$$\CA=\lim_{\longleftarrow}\CbarA{d}=\prod_{d\geq 0} C_d(\A)$$

For any $d\in\IN$, set
$$\CbarA{d}_+=\{f\in\CbarA{d}|\,f(0)=0\}
\aand
\CbarA{d}^\times=\{f\in\CbarA{d}|\,f(0)=1\}.$$
%and define similarly $\CA_+$ and $\CA^\times$.
Then, $\CbarA{d}_+$ is a Lie subalgebra of $\CbarA{d}$ and
$\CbarA{d}^\times$ is a subgroup of the group of invertible
elements in $\CbarA{d}$ since $f\in \CbarA{d}$ is invertible if,
and only if, $f(0)\neq 0$. The following is standard.

\begin{lemma}\label{le:exp log}
The standard exponential and logarithm functions
$$\exp(x)=\sum_{n\geq 0}\frac{x^n}{n!} \aand
\log (y)=\sum_{n\geq 1}\frac{(-1)^{n-1}}{n}(y-1)^n$$
yield well--defined maps
$$\exp:\CbarA{d}_+\to\CbarA{d}^\times\aand
\log:\CbarA{d}^\times\to\CbarA{d}_+$$
%$$\exp:\CA_+\to\CA^\times\aand
%\log:\CA^\times\to\CA_+$$
which are each other's inverse.
\end{lemma}

Similarly, one can define a Lie subalgebra $\CA_+$ of $\CA$ and
a subgroup $\CA^\times$ of the group of invertible elements of $\CA$.
By Lemma \ref{le:exp log}, the exponential and logarithm functions
give mutually inverse maps $\CA_+\rightleftarrows\CA^\times$.

\subsection{The pro--unipotent group $\NA$}
%------------------------------------------------------

The completion $\CA$ inherits a bialgebra structure from $C(\A)$ since
$$\Delta(C_d(\A))\subset\bigoplus_{a+b=d}C_a(\A)\otimes C_b(\A)$$
Set $\n_{>d}(\A)=\n(\A)\cap C_{>d}(\A)$ and
$$\nbarA{d}=\n(\A)/\n_{>d}(\A)$$
Then, the Lie algebra
$$\wh{n}(\A)=
\lim_{\longleftarrow} \nbarA{d}=\prod_{d\geq 0} \n_d(\A)$$
is the subspace of primitive elements in $\CA$ and a Lie subalgebra
of $\CA_+$ by Lemma \ref{le:primitive}.

Define the {\it Ringel--Hall group} $\wh{N}(\A)$ of $\A$ to be the
set of grouplike elements of $\CA$. Since any grouplike element
$f\in\CA$ satisfies $f(0)=\eta(f) =1$, this is a subset of the set
$\CA ^\times$ of invertible elements. It is easy to check using
the bialgebra property that it is also a subgroup, i.e. is closed
under multiplication.

\begin{prop}\hfill
\begin{enumerate}
\item The exponential and logarithm maps restrict to isomorphisms
$$\exp:\nA\to\NA\aand\log:\NA\to\nA$$
\item The group $\wh{N}(\A)$ is a pro--unipotent group with Lie algebra
$\wh{n}(\A)$.
\end{enumerate}
\end{prop}
\begin{pf}
(i) readily follows from Lemma \ref{le:exp log} and the fact that $\Delta:
\CA\to\CA\otimes\CA$ is an algebra homomorphism (see, \eg \cite
[Thm. 3.2]{Reutenauer}).

(ii) Let $\pi_{\leq d}$ be the projection $C(\A)\to\CbarA{d}$. As an
abstract subgroup of $\CA$, $\NA$ is the inverse limit of the groups
$\NbarA{d}=\pi_{\leq d}(\NA)$. By (i),
$$%\begin{equation}\label{eq:NbarA}
\NbarA{d}=
\pi_{\leq d}(\NA)=
\pi_{\leq d}(\exp(\nA))=
\exp(\nbarA{d})
$$%\end{equation}
Since any element $f\in\CbarA{d}_+\supset\nbarA{d}$ is nilpotent,
Lemma \ref{le:exp log} implies that
$$\NbarA{d}=\{f\in\CbarA{d}^\times|\log(f)\in\nbarA{d}\}$$
is a Zariski closed, unipotent subgroup of $\CbarA{d}^\times$ with
Lie algebra $\nbarA{d}$. The conclusion follows since the projection
maps $\NbarA{d}\to\NbarA{d'}$, $d\geq d'$ are clearly regular.
\end{pf}

\subsection{The pro--solvable group $\wh{B}(\A)$}
%--------------------------------------------------------------

Let $\HA=\Hom_{\IZ}(\KA,\IC^*)$ be the torus of characters of
$\KA$. $\HA$ acts by
bialgebra automorphisms on $C(\A)$ by
$$\sum_{\gamma\in K_{\geq 0}(\A)} X_\gamma\mapsto
\sum_{\gamma\in K_{\geq 0}(\A)}\zeta(\gamma) X_\gamma.$$
where $\zeta\in \HA$. This action extends to $\CA$ and leaves
$\NA$ invariant. By definition, $\BA$ is the semidirect product
\begin{equation}\label{eq:BA}
\BA=
\HA\ltimes\NA=
\lim_{\longleftarrow} \HA\ltimes\NbarA{d}
\end{equation}
We refer to $\BA$ as the {\it extended Ringel--Hall group}
of $\A$. $\BA$ is a pro--solvable, pro--algebraic group with
maximal torus $\HA$. Its Lie algebra is
$$\wh{\bb}(\A)=\h(\A)\ltimes\wh{\n}(\A)$$
where $\hA=\Hom(\KA,\IC)$ is the Lie algebra of $\HA$.

\section{Stability conditions and wall--crossing}\label{se:stability conditions}
%===================================

\subsection{Stability conditions}
%-------------------------------------

We shall define a {\it stability condition} on $\A$ to be a group homomorphism
\[Z\colon K(\A) \to \C\]
such that $Z(K_{>0}(\A))\subset\IH$, where $\IH\subset\C$ is the upper
half--plane\footnote{This does not quite agree with the definition given
in \cite{Bridgeland} where $Z$ is allowed to take $K_{>0}(\A)$ into the
half--closed half--plane. The difference will be of no importance in this
paper however.}. Let $\Stab(\A)$ denote the set of all stability conditions
on $\A$. Since the positive cone $K_{>0}(\A)$ is generated by the classes
of the simple modules $S_1,\ldots,S_N$ there is a bijection
\[\Stab(\A) \isom \IH^N\]
sending a stability condition $Z$ to the $N$--tuple $(Z(S_1),\ldots, Z(S_N))$.
We may therefore regard $\Stab(\A)$ as a complex manifold.

Let $Z$ be a stability condition on $\A$. Each nonzero module $M\in \A$
has a {\it phase}
$$\phi(M)=\frac{1}{\pi}\arg Z(M) \in (0,1).$$
A module $M$ is said to be {\it $Z$--semistable} if it is nonzero and if
\[0\neq A\subset M  \implies \phi(A)\leq \phi(M).\]

\subsection{Wall--crossing}
%--------------------------------

For any pair of classes $\beta,\gamma\in K_{>0}(\A)$ which are not
proportional over $\IQ$, consider the real codimension one submanifold
$W_{\beta,\gamma}$ of $\Stab(\A)$ given by
\[W_{\beta,\gamma}=\{Z\in \Stab(\A): Z(\beta)/Z(\gamma)\in\R_{>0}\}.\]
$W_{\beta,\gamma}$ is known as a \emph{wall}.

Fix a class $\alpha\in K_{>0}(\A)$ and consider the walls $W_{\beta,
\gamma}$ as $\beta,\gamma\in K_{>0}$ vary over the finitely many pairs
such that $\beta+\gamma=\alpha$. The connected components of the
complement of these walls in $\Stab(\A)$ are called \emph{chambers}. It is clear that
in each chamber the set of semistable objects of type $\alpha$ is constant.
However the set of semistable objects may change as one crosses a wall
from one chamber to a neighbouring one. We refer to this behaviour as
\emph{wall--crossing}.

\subsection{The functions $\delta_\gamma$} % delta functions
%------------------------------------------------------

The following result shows that semistability \wrt a given
$Z\in\Stab(\A)$ is an open condition. Recall that a family of
$R$-modules over a base variety $S$ is a vector bundle $M$ on $S$
together with a ring homomorphism $R\to \End_S(M)$.

\begin{lemma}\label{le:open}
Given a family of modules over a variety $S$, the subset of points of $S$
which correspond to $Z$--semistable modules is open.
\end{lemma}
\begin{pf}
Fix a class $\gamma\in K_{>0}(\A)$. It is immediate from the definitions that
a module $M$ of class $\gamma$ is $Z$--semistable if, and only if it is $\theta
$--semistable in the sense of King \cite{King}, where $\theta\colon K(\A)\to\R$
is given by
\[\theta(\beta)=-\Im (Z(\beta)/Z(\gamma)).\]
In turn, King shows that $\theta$-semistability coincides with GIT semistability
for the action of a reductive algebraic group on an affine variety \cite[Proposition 3.1]
{King}. It follows from this that semistability is an open condition.
\end{pf}

Given $Z\in \Stab(\A)$ and $\gamma\in K_{> 0}(\A)$ define $\delta_\gamma
\in\H_\gamma(\A)$ to be the characteristic function of $Z$--semistable modules
of type $\gamma\in K(\A)$.
\[\delta_{\gamma}(M)
 =\begin{cases}
 \,1 &\text{if $M$ is $Z$--semistable and $[M]=\gamma$,} \\
 \, 0 &\text{otherwise.}
\end{cases}
 \]
Lemma \ref{le:open} implies that $\delta_\gamma$ is a constructible function.
Clearly $\delta_\gamma$ depends on $Z\in \Stab(\A)$ in a discontinuous way
because of wall--crossing behaviour.

\subsection{Harder--Narasimhan filtrations and Reineke inversion} % HN
%---------------------------------------------------------------------------------

Since $\A$ is a finite length category, the Harder--Narasimhan property
holds (in this context see for example \cite[Prop. 2.4]{Bridgeland}): every
nonzero module $M$ has a unique filtration
\[0=M_0\subset M_1 \subset \cdots \subset M_n=M\]
where the successive factors $F_i=M_i/M_{i-1}$ are $Z$--semistable and
of strictly decreasing phases
\[\phi(F_1)>\cdots >\phi(F_n).\]
The following result is due to Reineke \cite{Reineke}.

\begin{thm}\label{th:Reineke}
For every stability condition $Z\in\Stab(\A)$, and every $\gamma\in K_{>0}
(\A)$, the following holds in $\H(\A)$.
\begin{align}%{gather}
\kappa_{\gamma}&=
\sum_{n\geq 1}\sum_{\stackrel{\gamma_1+\cdots+\gamma_n=
\gamma}{\phi(\gamma_1)>\cdots>\phi(\gamma_n)}}
\delta_{\gamma_1}*\cdots*\delta_{\gamma_n},
\label{reineke}
\\[1.1em]
\delta_{\gamma}&=
\sum_{n\geq 1}
\sum_{\stackrel{\gamma_1+\cdots+\gamma_n=\gamma}
{\phi(\gamma_1+\cdots+\gamma_i)>\phi(\gamma)}}
(-1)^{n-1}\kappa_{\gamma_1}*\cdots*\kappa_{\gamma_n},
\label{eq:reinekeinversion}
\end{align}%{gather}
where the (finite) sums are over elements $\gamma_i\in K_{>0}(\A)$.
\end{thm}
\begin{pf}
The first identity follows immediately from the definition of the product in
$\H(\A)$ and the existence and uniqueness of Harder--Narasimhan filtrations.
The second is proved in \cite[Theorem 5.1]{Reineke}
\end{pf}

The following is a direct consequence of \eqref{eq:reinekeinversion}
\begin{cor}
The elements $\delta_\gamma$ lie in the subalgebra $C(\A)\subset\H(\A)$
generated by the elements $\kappa_{\gamma}$.
\end{cor}

\subsection{The elements $\SS_\ell$}\label{ss:SS ell}
%---------------------------------------------

Let $Z\in \Stab(\A)$ be a stability condition. Given a ray $\ell=\IR_{>0}
e^{i\pi\phi}$, the infinite sum
\[\SS_\ell=1+
\sum_{\substack{\gamma\in K_{>0} (\A)\\ Z(\gamma)\in \ell}}\delta_\gamma\]
defines an element of $\CA$ which has the value 1 on a module $M$
if it is zero or semistable of phase $\phi$, and has value zero otherwise.

\begin{lemma}
The element $\SS_\ell$ is grouplike and therefore lies in $\NA$.
\end{lemma}
\begin{pf}
This amounts to the statement that a module $M\oplus N$ is semistable of phase
$\phi$ precisely if $M$ and $N$ are. This is a standard fact but for the reader's
convenience we sketch the proof. The only non--obvious implication is that if $M$
and $N$ are semistable of phase $\phi$ then so is $M\oplus N$. Suppose $A\subset
M\oplus N$ is a subobject with phase $\phi(A)>\phi$. We can assume that $A$ is
semistable: otherwise pass to a submodule of larger phase and repeat. The
inclusion $A\subset M\oplus N$ gives a nonzero map $A\to M$ or $A\to N$. But
the image of such a map  has phase larger than $\phi$ (because it is a quotient
of $A$ which is semistable of phase $>\phi$) and smaller than $\phi$ (because
it is a subobject of $M$ or $N$ which are semistable of phase $\phi$). This gives
a contradiction.
\end{pf}

\subsection{The element $\Splus$}
%------------------------------------------

Let now $\Splus\in\CA$ be the element given by
\[\Splus=\sum_{\gamma\in K_{\geq  0}(\A)} \kappa_\gamma.\]
The function $\Splus$ takes the value 1 on every module. It is
therefore grouplike and lies in $\NA$.

\subsection{Clockwise product}
%-------------------------------------

The \HN relation \eqref{reineke} can be written more compactly in terms
of the elements $\SS_\ell$ and $\Splus$ as follows
\begin{equation}\label{eq:Hall clockwise}
\Splus=\stackrel{\curvearrowright}{\prod_\ell}\SS_{\ell}
\end{equation}
where the product is an infinite product in the group $\wh{N}(\A)$ over all
rays of the form $\ell=\R_{>0}Z(\gamma)$ for some $\gamma\in K_{>0}(\A)$,
taken in clockwise order. This product makes sense in $\wh{N}(\A)$ because
it is finite when evaluated on modules $M$ of a fixed type $\gamma\in K_
{\geq 0}(\A)$.

\subsection{Joyce's elements $\mathbf{\epsilon_\alpha}$}
%-----------------------------------------------------------------------

Given a class $\alpha\in K_{>0}(\A)$, Joyce defines an element $\epsilon
_\alpha\in C_\alpha(\A)$ by the finite sum \cite{Joyce}
\begin{equation}
\label{joycelog}
\epsilon_\alpha=
\sum_{n\geq 1}
\sum_{\stackrel{\gamma_1+\cdots+\gamma_n=\alpha}{Z(\gamma_i)\in\R
_{>0}Z(\alpha)}}
\frac{(-1)^{n-1}}{n}\,\delta_{\gamma_1}*\cdots *\delta_{\gamma_n}
\end{equation}
He then shows that $\epsilon_\alpha$ is supported on indecomposable
modules and hence defines an element  $\epsilon_\alpha\in\n_\alpha(\A)$.
From our Hopf algebraic perspective, this is clear: just as in formula \eqref
{log}, the equation \eqref{joycelog} expresses the grouplike element $\SS
_\ell$ as an exponential
\begin{equation}\label{eq:joyceexp}
\SS_\ell=
\exp\bigg(\sum_{\substack{\alpha\in K_{>0}(\A)\\ Z(\alpha)\in \ell}}\epsilon
_\alpha\bigg).
\end{equation}
Thus, the elements $\epsilon_\alpha$ are primitive and, by Lemma \ref
{le:primitive} are supported on non--zero indecomposable objects.

\subsection{}
%--------------

The relations between the elements $\epsilon, \delta$ and $\kappa$ are summarised in the following diagram.
\[\xymatrix{(\epsilon)  \ar@/_1pc/[rr]_{\stackrel{\text{\scriptsize exp}}{\eqref{exp}}}
&& (\delta) \ar@/_1pc/[ll]_{\stackrel{\text{\scriptsize log}}{\eqref{joycelog}}} \ar@/_1pc/[rr]_{\stackrel
{\text{\scriptsize Harder-Narasimhan}}{\eqref{reineke}}}
&&(\kappa)\ar@/_1pc/[ll]_{\stackrel{\text{\scriptsize Reineke inversion}}{\eqref{eq:reinekeinversion}}}
  }\]
Note that it is identical to the right--hand part of the diagram of Section \ref{dia}.

\section{Stability conditions and Stokes data}\label{se:stability and Stokes}
%==================================

%This section contains the main result of this paper (Theorem \ref{th:le main}).
We show in this section that a stability condition $Z$ on $\A$ defines Stokes
data for an irregular connection on $\IP^1$ with values in the Ringel--Hall Lie
algebra of $\A$. We show moreover that this connection varies isomonodromically
with $Z$.

\subsection{Irregular $\BA$--connections on $\IP^1$}
%-----------------------------------------------------------------

We begin by adapting the definition of Stokes data to irregular connections
with values in the pro--solvable group $\BA$ defined in Section \ref{se:Hall}.

Let
$$\Phi=\{\alpha\in K_{>0}(\A)|\n_\alpha(\A)\neq 0\}\subset\hA^*$$
be the set of roots of $\wh{B}(\A)$ relative to the torus $H(\A)$ and set
$$\h(\A)_\reg=\h(\A)\setminus\bigcup_{\alpha\in\Phi}\Ker(\alpha)$$

Let $\wh{P}$ be the holomorphically trivial, principal $\BA$--bundle on
$\IP^1$. By this we mean the following: the group $\BA$ is the inverse
limit of the solvable algebraic groups $\BbarA{d}$ and $\wh{P}$ is the
limit of the corresponding principal bundles $\PbarA{d}$. In particular,
a section of $P$ is holomorphic if the induced section of each $\PbarA
{d}$ is.

Fix $Z\in\h(A)_\reg$ and consider a connection on $\wh{P}$ of the form
$$\nabla=d-\left(\frac{Z}{t^2}+\frac{f}{t}\right)dt.$$
where $f\in\nA$. $\nabla$ is the inverse limit of the connections
$$\nabla_{\leq d}=d-\left(\frac{Z}{t^2}+\frac{\pi_{\leq d}(f)}{t}\right)dt.$$
on $\PbarA{d}$, where $\pi_{\leq d}:\wh{\n}(\A)\to\nbarA{d}$ is the projection,
and each $\nabla_{\leq d}$ satisfies the assumptions of Section \ref
{ss:assumptions}.
Indeed, $\HA$ is a maximal torus in $\BbarA{d}$ with corresponding
root space decomposition \eqref{eq:root decomposition 1} given by
$$\bbarA{d}=
\h(A)%\oplus
\bigoplus_{\substack{\alpha\in\Phi\\\dim\alpha\leq d}}\n_{\alpha}(\A),
$$
$Z\in\hA$ is a regular element and the projection of $\pi_{\leq d}(f)$
onto $\hA$ is zero.

\subsection{Canonical fundamental solutions}
%-------------------------------------------------------

As in Section \ref{ss:Stokes factors}, we define a Stokes ray of $\nabla$
to be a ray of the form $\ell=\R_{>0}Z(\alpha)$ for $\alpha\in\Phi$. Note
that there may be an infinite number of such rays. A ray is admissible if
it is not a Stokes ray. Since $\n(\A)$ is isomorphic to $\nbarA{d}\oplus\n_{>d}
(\A)$ as $\HA$--module, an admissible ray for $\nabla$ is admissible for
each $\nabla_{\leq d}$. One can therefore use the existence and uniqueness
statement of Theorem \ref{jurk} to deduce the following
\begin{thm}
Given an admissible ray $r$, there is a unique holomorphic fundamental
solution
\[Y_r\colon \IH_r\to \BA\]
of $\nabla$ such that $Y_r(t)\cdot\exp(Z/t)\to 1$ as $t\to 0$ in $\IH_r$.
\end{thm}

\subsection{Stokes factors}
%--------------------------------

The definition of the Stokes factors of $\nabla$ requires a little care
since the set of Stokes rays of $\nabla$ need not be discrete. If $r_1
\neq -r_2$ are two admissible rays however, ordered so that the closed
sector $\ol{\Sigma}(r_1,r_2)\subset\IC^*$ swept by clockwise rotation
from $r_1$ to $r_2$ is convex, there is a unique element $S_{\ol{\Sigma}
(r_1,r_2)}\in\BA$ such that
$$Y_{r_2}=Y_{r_1}\cdot S_{\ol{\Sigma}(r_1,r_2)}$$
on $\IH_{r_1}\cap \IH_{r_2}$. By Proposition \ref{pr:spectral},
$S_{\ol{\Sigma}(r_1,r_2)}$ is of the form $\exp(X)$ where
$$X\in
\prod_{\substack{\alpha\in \Phi\\[.3 ex] Z(\alpha)\in\ol{\Sigma}(r_1,r_2)}}
\n_\alpha(\A)$$

\begin{defn}
$\nabla$ admits a Stokes factor $S_\ell\in\BA$ along the Stokes ray $\ell$
if the elements $S_{\ol{\Sigma}(r_1,r_2)}$ tend to $S_\ell$ as the admissible
rays $r_1,r_2$ tend to $\ell$ in such a way that $\ell\in\ol{\Sigma}(r_1,r_2)$.
\end{defn}

\begin{prop}\hfill
\begin{enumerate}
\item The connection $\nabla$ admits a Stokes factor $S_\ell$ along any
Stokes ray $\ell$.
\item Given two admissible rays $r_1,r_2$ as above, one has
$$S_{\ol{\Sigma}(r_1,r_2)}=\clockwise_{\ell\subset\ol{\Sigma}(r_1,r_2)}S_\ell$$
\end{enumerate}
\end{prop}
\begin{pf}
Both statements hold for each solvable quotient $\BbarA{d}$ of
$\BA$.
\end{pf}

\subsection{Stokes multipliers}
%------------------------------------

Unlike the definition of the Stokes factors, that of the Stokes multipliers
$S_\pm$ of $\nabla$ relative to the choice of a ray $r$ such that $\pm
r$ are admissible is straightforward and given, as in Section \ref{ss:multipliers}
by
$$Y_{r,\pm}(t)=Y_{-r}(t)\cdot S_\pm,\quad t\in\IH_{-r}$$
where $Y_{r,+}$ and $Y_{r,-}$ are the analytic continuations of $Y_r$
to $\IH_{-r}$ in the anticlockwise and clockwise directions respectively.
By Lemma \ref{veryeasy}, $S_\pm$ are given by the clockwise
products over Stokes factors
\begin{equation}\label{eq:another clockwise}
S_+=\clockwise_{\ell\in\IH_{ir}}S_\ell\aand
(S_-)^{-1}=\clockwise_{\ell\in\IH_{-ir}}S_\ell
\end{equation}

\subsection{Stability conditions}
%-------------------------------------

Assume now that $Z\in\Stab(\A)\subset\hA$ is a stability condition.
Note that $Z\in\hA_\reg$ since $Z(\alpha)\in\IH$ for any $\alpha\in
K_{>0}(\A)$.

For any ray $\ell=\IR_{>0}e^{i\pi\phi}$, let $\SS_\ell$ be the
characteristic function of semistable objects of phase $\phi$
defined in Section \ref{ss:SS ell}. By \eqref{eq:joyceexp},
$$\SS_\ell=
\exp\left(\sum_{\alpha:Z(\alpha)\in\ell}\epsilon_\alpha\right)$$
where $\epsilon_\alpha\in\n(\A)_\alpha$. This shows in particular
the following

\begin{lemma}\label{le:reassuring}
If $\ell$ is not a Stokes ray of $\nabla$, then $\SS_\ell=1$.
\end{lemma}

Thus, only the elements $\SS_\ell$ corresponding to Stokes
rays of $\nabla$ are non--trivial.

\subsection{}
%--------------

The following is the main result of this paper.

\begin{thm}\label{th:le main}\hfill
\begin{enumerate}
\item
There exists a unique connection $\nabla_{\A,Z}$ of the form
\begin{equation}\label{eq:Hall nabla}
\nabla_{\A,Z}=d-\left(\frac{Z}{t^2}+\frac{f}{t}\right)dt.
\end{equation}
whose Stokes data is given by either of the following equivalent
conditions:
\begin{enumerate}
\item The Stokes factor corresponding to a Stokes ray $\ell=
\IR_{>0}\exp(i\pi\phi)$ is the characteristic function $\SS_\ell$
of $Z$--semistable modules of phase $\phi$.
\item The Stokes multipliers $S_+,S_-$ relative to the ray $r=
\IR_{>0}$ are the function $\Splus$ which takes the value 1
on every module, and the identity element $1_0$ respectively.
\end{enumerate}
\item
The components of $f=\sum_{\alpha\in K_{>0}(\A)}f_\alpha$ are given
by the Lie series
$$f_\alpha=
\sum_{n\geq 1}
\sum_{\substack{\alpha_i\in K_{>0}(\A)\\ %\alpha_1,\ldots,\alpha_n\in K_{>0}(\A)\\
\alpha_1+\cdots+\alpha_n=\alpha}}
J_n(Z(\alpha_1),\ldots, Z(\alpha_n))\,
\epsilon_{\alpha_1}*\cdots*\epsilon_{\alpha_n}.$$
where the $J_n$ are the functions appearing in Theorem
\ref{th:Stokes inverse}.

\item
As $Z$ varies in $\Stab(\A)$, the family of connections $\nabla_{\A,Z}$
varies isomonodromically. In particular, $f_\alpha(Z)$ is a holomorphic
function of $Z$ and satisfies the PDE
$$df_\alpha=
\sum_{\beta+\gamma=\alpha}[f_\beta,f_\gamma]d\log\gamma$$
\end{enumerate}
\end{thm}
\begin{pf}

(i) We first show the equivalence of (a) and (b). Note that since $Z\in\Stab
(\A)$, the Stokes rays of a connection $\nabla$ of the form \eqref{eq:Hall nabla}
are contained in the upper half plane $\IH$. In particular the ray $r=\IR_{>0}$
is such that $\pm r$ are admissible and the corresponding Stokes multiplier
$S_-$ is trivial by \eqref{eq:another clockwise}.

Assume that (a) holds and let $\calR$ be the set of Stokes rays of
$\nabla$. By \eqref{eq:another clockwise}
$$S_+
=
\clockwise_{\ell\in\calR:\ell\subset\IH}S_\ell
=
\clockwise_{\ell\in\calR:\ell\subset\IH}\SS_\ell
=
\clockwise_{\ell\subset\IH}\SS_\ell
=
\Splus$$
where the third identity follows from Lemma \ref{le:reassuring}
and the last one from the \HN relation \eqref{eq:Hall clockwise}.

% (ii)=>(i)
The implication (b)$\Rightarrow$(a) follows from \eqref
{eq:another clockwise}, \eqref{eq:Hall clockwise} and the fact
that the relations \eqref{reineke} can be inverted.

The existence and uniqueness of $f$, and the fact that it is given
by the Lie series (ii) follows from Theorem \ref{th:Stokes inverse}
(note that the series \eqref{main2} is finite since any $\alpha\in
K_{>0}(\A)$ can only be decomposed as the sum of elements
of $K_{>0}(\A)$ in finitely many ways. Thus, the Taylor series
of $\calS^{-1}$ yields a global inverse of the Stokes map in this
case).

%IMD
(iii)
The first assertion follows from the fact that, by condition (b), the
Stokes multipliers $S_\pm$ are constant functions of $Z$. The
second follows from Theorem \ref{th:IMD}. Indeed, since the
Stokes map has a global inverse, $f=f(Z)$ varies holomorphically
in $Z$ by Remark \ref{rk:IMD} and satisfies the isomonodromy
equations \eqref {iso}.
\end{pf}

\end{document}